\documentclass{article}
\usepackage{amsmath}
\usepackage{amsfonts}

\advance \textheight by \topskip
\advance \textheight by 1pt
\makeatother
\def\bea{\begin{eqnarray}}
\def\ena{\end{eqnarray}}
\def\non{\nonumber}

\def\lar{\longrightarrow}
\def\ra{\rightarrow}
\def\sea{\searrow}
\def\deg{\hbox{deg}}
\def\rk{\hbox{rank}}
\def\Ker{\hbox{Ker}}
\def\tr{\hbox{tr}}
\def\Im{\hbox{Im}}
\def\ch{\hbox{ch}}
\def\sgn{\hbox{sgn}}
\def\Res{\text{Res}}

\def\lar{\longrightarrow}
\newtheorem{prop}{Proposition}
\newtheorem{theorem}{Theorem}

\newtheorem{lemma}{Lemma}
\newtheorem{cor}{Corollary}

\def\pf{{\it Proof.} \quad}
\newcommand{\bc}[2]{
\left(
\begin{array}{c}{#1}\\{#2}\end{array}
\right)}
\newcommand{\qbc}[2]{
\left[
\begin{array}{c}{#1}\\{#2}\end{array}
\right]_q}
\newcommand{\pbc}[2]{
\left[
\begin{array}{c}{#1}\\{#2}\end{array}
\right]_p}
\newcommand{\qqbc}[2]{
\left[
\begin{array}{c}{#1}\\{#2}\end{array}
\right]_{q^2}}
\newcommand{\tetc}[4]{
\left(
\begin{array}{ccc}
{}&{#1}&{}\\
{#2}&{#3}&{#4}
\end{array}
\right)}
\newcommand{\qqtetc}[4]{
\left[
\begin{array}{ccc}
{}&{#1}&{}\\
{#2}&{#3}&{#4}
\end{array}
\right]_{q^2}}

\title{
Residues of $q$-Hypergeometric Integrals
and Characters of Affine Lie Algebras
}

\author{
Atsushi Nakayashiki\thanks{
Faculty of Mathematics,
Kyushu University,
Ropponmatsu 4-2-1, Fukuoka 810-8560, Japan, \quad 
e-mail: 6vertex@math.kyushu-u.ac.jp}
}

\date{}
\begin{document}
\maketitle
\begin{abstract}
We study certain subspaces of solutions to the $sl_2$ rational qKZ equation 
at level zero.
Each subspace is specified by the vanishing of the residue at a certain divisor
which stems from models in two dimensional integrable field theories.
We determine the character of the subspace which is parametrized by 
the number of variables and the $sl_2$ weight of the solutions. 
The sum of all characters with a fixed weight gives rise to the branching 
functions of the irreducible representations of $sl_2$ in the 
level one integrable highest weight representations of $\widehat{sl_2}$. 
It is written in the fermionic form.
\end{abstract}

\section{Introduction}
In this paper we study a problem arising from the study
of two dimensional integrable quantum field theories. 
It is related with the calculation of the character 
of the space of local operators.

Let us first recall the known simplest example to give a preliminary 
insight on what kind of problem we are going to study.
For the two dimensional massive Ising field theory in zero magnetic field, the problem of classifying local operators reduces to classifying 
the solutions to the equation \cite{CM,Ch,K}
\bea
&&
f(x_1,\cdots,x_{n-2},x,-x)=0,
\label{eq0-0}
\ena
where $f(x_1,\cdots,x_n)$ is a symmetric polynomial.
If we denote $R_n$ the ring of symmetric polynomials of $x_j$'s, 
the space of solutions of (\ref{eq0-0}) is given by
\bea
&&
M_n:=R_n\Delta^{+}_n,
\quad
\Delta^{+}_n=\prod_{i<j}^n(x_i+x_j).
\non
\ena
Specifying $\deg\, x_i=1$, $\forall i$, $M_n$ decomposes into the direct sum 
of the degree $d$ subspaces $M_n(d)$.
The generating function of the dimension of $M_n(d)$ is calculated as
\bea
&&
\ch\, M_n:=\sum_d q^d \dim M_n(d)=\frac{q^{\frac{1}{2}n(n-1)}}{(q)_n},
\quad
(q)_n=\prod_{j=1}^n(1-q^j).
\non
\non
\ena
Let $M^{(i)}$, $i=0,1$ be the direct sum of $M_n$ with $n$ being even for $i=0$
and odd for $i=1$.
Then
\bea
&&
\ch\, M^{(i)}=\sum_{n\equiv i\, mod.\, 2}\frac{q^{\frac{1}{2}n(n-1)}}{(q)_n}.
\non
\ena
These infinite sums can be rewritten into the infinite product form \cite{Ch,K}
\bea
&&
\ch\, M^{(0)}=\ch\, M^{(1)}=\prod_{j=1}^\infty (1+q^j).
\non
\ena
The last expression is the character of the Fock space of free fermions.
The result is expected because the two dimensional Ising model in a zero
magnetic field is described, at the critical point, by the free fermion 
conformal field theory.

The two particle $S$-matrix of the Ising model is $-1$. 
To go beyond the Ising model there are two possibilities.
One is to consider models with more complicated but diagonal $S$-matrices.
In this case similar calculations were carried for several models in \cite{K}. 
The other is to consider models with non-diagonal $S$-matrices.
Up to now no similar calculations have been successful in such cases.

In this paper we carry out similar calculations for the simplest model with 
a non-diagonal $S$-matrix, the $SU(2)$ invariant Thirring model \cite{Smir2}. 
In this case the symmetric polynomial $f$ in (\ref{eq0-0}) is replaced by
a tensor valued function satisfying certain system of equations. 
The solutions to this system of equations are given by 
multi-dimensional $q$-hypergeometric integrals. 
For these integrals we have to solve the equation 
corresponding to (\ref{eq0-0}).
Thus the problem we consider looks far more complicated than the Ising case.
Nevertheless we can solve this problem completely as shown in this paper.
Let us explain our results supplying more precise statements.

Let $V$ be a two dimensional vector space and $S(\beta)$ the linear operator
acting on $V^{\otimes 2}$ defined by
\bea
&&
S(\beta)=S_0(\beta)\hat{S}(\beta),
\quad
\hat{S}(\beta)=
\frac{\beta-\pi i P}{\beta-\pi i},
\quad
S_0(\beta)=
\frac{\Gamma(\frac{\pi i+\beta}{2\pi i})\Gamma(\frac{-\beta}{2\pi i})}
{\Gamma(\frac{\pi i-\beta}{2\pi i})\Gamma(\frac{\beta}{2\pi i})},
\non
\ena
where $P$ is the operator which permutes the components of $V^{\otimes 2}$. 
In general for a linear operator $A$ acting on $V^{\otimes 2}$, 
we denote $A_{ij}$ the operator on $V^{\otimes n}$
which acts $i$-th and $j$-th components as $A$ and 
other components as an identity.

Consider the system of equations for a $V^{\otimes n}$-valued function $f$:
\bea
&&
P_{i,i+1}S_{i,i+1}(\beta_i-\beta_{i+1})
f(\beta_1,\cdots,\beta_n)
=f(\cdots,\beta_{i+1},\beta_i,\cdots),
\label{eq0-1}
\\
&&
P_{n-1,n}P_{n-1,n-2}\cdots P_{1,2}
f(\beta_1-2\pi i,\cdots,\beta_n)
=
(-1)^{\frac{n}{2}}f(\beta_2,\cdots,\beta_n,\beta_1).
\label{eq0-2}
\ena
This set of equations implies the $sl_2$ rational qKZ equation at 
level zero \cite{FR,Smir0} and it can be solved in terms of multi-dimensional
$q$-hypergeometric integrals \cite{KS,NPT,Smir1,Smir2}. 
In this description of solutions,
meromorphic solutions of (\ref{eq0-1}) and (\ref{eq0-2}) are parametrized
by anti-symmetric polynomials $P(X_1,\cdots,X_\ell)$ of $X_a$'s of degree
less than $n$ in each variable with the
coefficients in the space $C_n$ of symmetric and $2\pi i$-periodic meromorphic 
functions of $\beta_j$'s, here $\ell$ specifies the $sl_2$ weight of 
a solution. We restrict ourselves to the case where 
the coefficients are in the space of symmetric polynomials of 
$x_j=e^{\beta_j}$, $1\leq j\leq n$ and write $P$ as 
$P(X_1,\cdots,X_\ell|x_1,\cdots,x_n)$.
We denote by $\Psi_P$ the solution corresponding to
$P$. It is known \cite{NT} that $\Psi_P$, 
as a function of $\beta_n$,
has at most a simple pole at $\beta_n=\beta_{n-1}+\pi i$.
Then the condition corresponding to (\ref{eq0-0}) is
\bea
&&
\Res_{\beta_n=\beta_{n-1}+\pi i}\, f(\beta_1,\cdots,\beta_n)=0.
\label{eq0-3}
\ena
This condition should be rewritten as the condition for $P$.
In fact we have derived equations for $P$ which imply (\ref{eq0-3}).
They are
\bea
&&
P(X_1,\cdots,X_{\ell-1},\pm x^{-1}|x_1,\cdots,x_{n-2},x,-x)=0.
\label{eq0-4}
\ena
We conjecture that these equations are equivalent to (\ref{eq0-3}).

Let $U_{n,\ell}$ be the space of solutions of (\ref{eq0-4}).
Since $P$ is an anti-symmetric polynomials of $X_a$'s, it can be considered
as an element of the $\ell$-th exterior product of a free $R_n$-module with
the monomials $X^j$ as a basis:
\bea
&&
P\in \wedge^\ell H^{(n)},
\quad
H^{(n)}=\oplus_{j=0}^{n-1} R_n X^j.
\non
\ena
There exist two special solutions 
\bea
&&
\Xi_1\in U_{n,1},
\quad
\Xi_2\in U_{n,2},
\non
\ena
for which $\Psi_{\Xi_1}$ and $\Psi_{\Xi_2}$ vanish identically 
\cite{NPT,Smir3,Tar}.
Moreover a polynomial $P$ with the coefficients in $C_{n}$ such that 
$\Psi_P$ vanishes identically belongs to the space \cite{Tar}
\bea
&&
\wedge^{\ell-1} \hat{H}^{(n)}\wedge \Xi_1
+
\wedge^{\ell-2} \hat{H}^{(n)}\wedge \Xi_2,
\quad
\hat{H}^{(n)}=\oplus_{j=0}^{n-1}C_{n} X^j.
\non
\ena
With this respect we consider the space
\bea
&&
M_{n,\ell}:=
\frac{U_{n,\ell}}
{U_{n,\ell-1} \wedge \Xi_1+U_{n,\ell-2}\wedge \Xi_2},
\non
\ena
and, for a non-negative integer $\lambda$,
\bea
&&
M^{(i)}_{\lambda}=
\oplus_{n\equiv i\, mod.\, 2,\, n-2\ell=\lambda} M_{n,\ell}.
\non
\ena
We define a degree, denoted by $\deg_1$, assigning 
$\deg_1\, X_{a}=-1$, $\deg_1\, x_j=1$. 
If $P$ is homogeneous by $\deg_1$, $\Psi_P$ becomes homogeneous of degree
\bea
&&
\deg_2\, P=\frac{n^2}{4}+\deg_1\, P,
\quad
P\in U_{n,\ell}.
\non
\ena
We introduce a grading on $M^{(i)}_{\lambda}$ by $\deg_2$. Then

\begin{theorem}
\bea
&&
\ch\, M^{(i)}_{\lambda}=
\sum_{n-2\ell=\lambda,\, n\equiv i\, mod.\, 2}
\frac{q^{\frac{n^2}{4}}}{(q)_n}
\Big(
\qbc{n}{\ell}-\qbc{n}{\ell-1}
\Big),
\label{branching}
\ena
where
\bea
&&
\qbc{n}{\ell}=\frac{(q)_n}{(q)_\ell(q)_{n-\ell}}.
\non
\ena
\end{theorem}

For a non-negative half integer $S$ let
\bea
&&
\chi_S(z)=\frac{z^{2S+1}-z^{-(2S+1)}}{z-z^{-1}}
\non
\ena
be the character of the $2S+1$-dimensional irreducible 
representation of $sl_2$. An easy calculation shows 
\bea
&&
\sum_{\lambda\equiv i\, mod.\, 2}
\chi_{\lambda/2}\, \ch\, M^{(i)}_{\lambda}
=
\sum_{n\equiv i\, mod.\, 2}\sum_{\ell=0}^n
\frac{z^{n-2\ell}q^{\frac{n^2}{4}}}{(q)_\ell\, (q)_{n-\ell}}.
\non
\ena
The last expression is the fermionic form of the character of 
the level one integrable highest weight representation $V(\Lambda_i)$ 
of the affine Lie algebra $\widehat{sl_2}$ presented in \cite{KMM,Me}. 
The formula (\ref{branching}) gives the branching function of the
$\lambda+1$-dimensional irreducible representation of $sl_2$ in 
$V(\Lambda_i)$. 

Finally we remark that, for the Ising model, it is possible to construct
a local operator, in the form of a set of form factors, to each element 
of $M^{(i)}$ \cite{CM,Ch}. 
The problem of constructing form factors corresponding to the elements
of $M^{(0)}_\lambda$ will be studied in the subsequent paper.

The present paper is organized in the following manner.
In section 2 the null residue equation (\ref{eq0-4}) for a polynomial $P$ 
is derived. From section 3 to section 9 is devoted to the case $n$ even.
We study the solution of (\ref{eq0-4}) with $\ell=1$ in section 3.
The module $U_{2n,1}$ is proved to be a free $R_{2n}$-module by giving 
its basis explicitly.
In section 4 the space $U_{2n,2}$ is studied. Again we prove that $U_{2n,2}$ is a free $R_{2n}$-module by constructing a basis.
In section 5 we present a candidate of a basis of $U_{n,\ell}$ in the space of 
exterior product of solutions for $\ell=1$ and $2$.
Some combinatorial identity related with $q$-binomials and $q$-tetranomials
are proved as a preparation of the study of $U_{n,\ell}$ for general $\ell$. 
In section 6 the structure theorem of $U_{2n,\ell}$ is given and a proof of it
is presented assuming the linear independence of the basis.
Sections 7 and 8 are devoted to the proof of the linear 
independence of the basis. 
The character formula (\ref{branching}) with $n$ even is derived in section 9.
In section 10 the results on the case $n$ odd are presented and proved.
In appendix A the details of calculations of the residue in section 2 is given.

\section{Null residue equation}
Let $R_{n}$ be the ring of symmetric polynomials of $x_1,...,x_n$.
We consider the space of polynomials in $X$ of degree at most $n-1$
with the coefficients in $R_{n}$:
\bea
&&
H^{(n)}:=\oplus_{k=0}^{n-1} R_{n} X^k.
\non
\ena
The $\ell$-th exterior product $\wedge^\ell H^{(n)}$ is 
identified with the space of anti-symmetric
polynomials $P$ in $X_1$, ..., $X_\ell$ of degree less than $n$ 
in each variable with the coefficients in $R_n$.
The identification is given by the map
\bea
X^{i_1}\wedge\cdots\wedge X^{i_\ell}
\mapsto
\text{Asym}(X_1^{i_1}\cdots X_\ell^{i_\ell})
=\sum_{\sigma\in S_\ell}\text{sgn}\,\sigma\, X_{\sigma(1)}^{i_1}\cdots 
X_{\sigma(\ell)}^{i_\ell}.
\non
\ena
In the following we suppose that the variables $X_a$, $x_j$ and
$\alpha_a$, $\beta_j$ are related by $X_a=e^{-\alpha_a}$, $x_j=e^{\beta_j}$.
To each element $P\in \wedge^\ell H^{(n)}$ and the index set
$M=(m_1,\cdots,m_\ell)$, $1\leq m_1<\cdots<m_\ell\leq n$, we associate 
the integral:
\bea
&&
I_M(P):=
\int_{C^\ell}\prod_{a=1}^\ell d\alpha_a \prod_{a=1}^\ell\phi_n(\alpha_a)
g_M\frac{P(X_1,\cdots,X_\ell|x_1,\cdots,x_n)}
{\prod_{a=1}^\ell\prod_{j=1}^n(1-X_ax_j)},
\label{eq2-1}
\ena
where 
\bea
&&
\phi_n(\alpha)=\prod_{j=1}^n
\frac{\Gamma(\frac{\alpha-\beta_j+\pi i}{-2\pi i})}
{\Gamma(\frac{\alpha-\beta_j}{-2\pi i})},
\quad
g_M=\prod_{a=1}^\ell
\Big(
\frac{1}{\alpha_a-\beta_{m_a}}
\prod_{j=1}^{m_a-1}
\frac{\alpha_a-\beta_j+\pi i}{\alpha_a-\beta_j}
\Big)
\prod_{1\leq a<b\leq \ell}(\alpha_a-\alpha_b+\pi i).
\non
\ena
The integration contour $C$ is defined as follows.
It goes from $-\infty$ to $+\infty$, separating two sets 
$\{\beta_j, \beta_j-2\pi i,\beta_j-4\pi i,\cdots|1\leq j\leq n\}$ and
$\{\beta_j-\pi i, \beta_j+\pi i,\beta_j+3\pi i,\cdots|1\leq j\leq n\}$.

The set of integrals $\{I_M(P)\}$ satisfies a certain system of equations.
Let us describe it. To this end we define a tensor valued function whose
component is $I_M(P)$.
Let $V$ be the vector space with the basis $v_{\pm}$:
\bea
&&
V=\mathbb{C}v_{+}\oplus \mathbb{C} v_{-}.
\non
\ena
To each $M=(m_1,\cdots,m_\ell)$, $1\leq m_1<\cdots<m_\ell\leq n$,
we define the element $v_M$ in $V^{\otimes n}$ by
\bea
&&
v_M=v_{\epsilon_1}\otimes\cdots\otimes v_{\epsilon_n},
\quad
M=\{j\,|\, \epsilon_j=-\,\}.
\non
\ena
Using them define the $V^{\otimes n}$-valued function $\psi_P$ by
\bea
&&
\psi_P=\sum_M I_M(P)v_M.
\non
\ena
To write down the equations we need the operator acting on 
$V^{\otimes 2}$ defined by
\bea
&&
\hat{S}(\beta)=\frac{\beta-\pi i P}{\beta-\pi i},
\quad
P(v_{\epsilon_1}\otimes v_{\epsilon_2})=
v_{\epsilon_2}\otimes v_{\epsilon_1}.
\non
\ena
In general for a linear operator $A$ acting on $V^{\otimes 2}$, writing
\bea
&&
A=\sum_a B_a\otimes C_a,
\quad
B_a,C_a\in \text{End}(V),
\non
\ena
we define the linear operator $A_{ij}$ acting on $V^{\otimes n}$ by
\bea
&&
A_{ij}=\sum_a\, 
1\otimes\cdots\otimes B_a\otimes\cdots \otimes C_a\otimes \cdots \otimes1,
\non
\ena
where
$B_a$ and $C_a$ are in i-th and j-th components respectively.

Then $\psi_P$ satisfies the following system of equations:
\bea
&&
P_{i,i+1}\hat{S}_{i,i+1}(\beta_i-\beta_{i+1})
\psi_P(\beta_1,\cdots,\beta_n)
=\psi_P(\cdots,\beta_{i+1},\beta_i,\cdots),
\label{eq2-2}
\\
&&
P_{n-1,n}P_{n-1,n-2}\cdots P_{1,2}
\psi_P(\beta_1-2\pi i,\cdots,\beta_n)
=
\psi_P(\beta_2,\cdots,\beta_n,\beta_1).
\label{eq2-3}
\ena

Here we remark on the solutions to the equations 
(\ref{eq0-1}) and (\ref{eq0-2}) in the introduction. 
By multiplying a scalar function we define 
$\Psi_P$ by
\bea
&&
\Psi_P=e^{\frac{n}{4}\sum_{j=1}^n\beta_j}
\prod_{j<k}\zeta(\beta_j-\beta_k)\psi_P,
\label{Psi-P}
\ena
where $\zeta(\beta)$ is a certain meromorphic function. 
See \cite{NT} and references therein for the precise description  
of $\zeta(\beta)$. We simply remark here that $\zeta(\beta)$ is holomorphic
and non zero at $\beta=-\pi i$.
The function $\Psi_P$ solves the equations (\ref{eq0-1}) and (\ref{eq0-2}).
Conversely every solution of (\ref{eq0-1}) and (\ref{eq0-2}) can be written
using $\Psi_P$. Let us give more detailed and precise statements on this fact.
Let $\Theta^{(n)}_{\pm}(X)=\prod_{j=1}^n(1\pm Xx_j)$ and set 
\cite{BBS,Smir3,Tar}
\bea
&&
2\Xi_1^{(n)}(X)=\Theta^{(n)}_{+}(X)+(-1)^{n-1}\Theta^{(n)}_{-}(X),
\label{Xi-1}
\\
&&
2\Xi_2^{(n)}(X_1,X_2)=
\Big(\Theta^{(n)}_{+}(X_1)\Theta^{(n)}_{+}(X_2)
-
\Theta^{(n)}_{-}(X_1)\Theta^{(n)}_{-}(X_2)
\Big)
\frac{X_1-X_2}{X_1+X_2}
\non
\\
&&
\qquad\qquad\qquad
+(-1)^n
\Big(
\Theta^{(n)}_{+}(X_1)\Theta^{(n)}_{-}(X_2)
-
\Theta^{(n)}_{+}(X_2)\Theta^{(n)}_{-}(X_1)
\Big).
\label{Xi-2}
\ena
Notice that the space of meromorphic solutions of (\ref{eq2-2}) 
and (\ref{eq2-3})
is a vector space over the field $C_n$ of symmetric and $2\pi i$-periodic
meromorphic functions of $\beta_1$,...,$\beta_n$. The ring $R_n$ becomes
a subring of $C_n$.
Consider the scalar extension of $H^{(n)}$ from $R_n$ to $C_n$:
\bea
&&
\hat{H}^{(n)}=C_{n}\otimes_{R_n} H^{(n)}=\oplus_{j=0}^{n-1}C_n X^j.
\non
\ena
Then

\begin{theorem}(\cite{Tar})\label{tar}
\begin{itemize}
\item[(i).] Any meromorphic solution of (\ref{eq0-1}) and (\ref{eq0-2}) which
takes the value in the space of $sl_2$ highest weight vectors in $V^{\otimes n}$ 
with the weight $n-2\ell$ can
be written as $\Psi_P$ for some $P\in \wedge^\ell \hat{H}^{(n)}$.

\item[(ii).] For $P\in \wedge^\ell \hat{H}^{(n)}$, $\Psi_P=0$ if and only if
\bea
&&
P\in 
\wedge^{\ell-1} \hat{H}^{(n)}\wedge \Xi_1^{(n)}
+
\wedge^{\ell-2} \hat{H}^{(n)}\wedge \Xi_2^{(n)}.
\non
\ena
\end{itemize}
\end{theorem}

The function $\Psi_P$ has a nice homogeneity property.
Let us define degrees, denoted by $\deg_1$, of elements in $\wedge^\ell H^{(n)}$by assigning 
\bea
&&
\deg_1\, X=-1,
\quad
\deg_1 
\, x_j=1.
\non
\ena
Then

\begin{prop}\label{prop2-1}
If $P$ is homogeneous by $\deg_1$, 
\bea
&&
\Psi_P(\beta_1+\theta,\cdots,\beta_n+\theta)
=e^{\deg_2(P)\theta}
\Psi_P(\beta_1,\cdots,\beta_n),
\non
\ena
where
\bea
&&
\deg_2(P)=\frac{n^2}{4}+\deg_1\, P.
\non
\ena
\end{prop}
\vskip2mm
\noindent
\pf Notice that $\alpha_a$, $\beta_j$ appear in the integrand of $I_M(P)$
in the form $\alpha_a-\beta_j$ or $\alpha_a-\alpha_b$ except in $P$.
The proposition immediately follows from this.
Q.E.D.
\vskip2mm

In the following the degree of an element of $\wedge^\ell H^{(n)}$ is measured by $\deg_1$ unless otherwise stated.

As a function of $\beta_n$, $\psi_P$ has at most a simple pole at 
$\beta_n=\beta_{n-1}+\pi i$ \cite{NT}.
We shall study the space of $P$ such that
\bea
&&
\Res_{\beta_n=\beta_{n-1}+\pi i}\psi_P=0.
\label{eq2-4}
\ena
By the remark below (\ref{Psi-P}) this is equivalent to (\ref{eq0-3}).
Let us find the condition for $P$ such that (\ref{eq2-4}) holds.

For a polynomial $P(X_1,\cdots,X_\ell|x_1,\cdots,x_{2n})$ 
of $x_i's$ and $X_j's$, define
\bea
&&
\rho_{\pm}(P)=P(X_1,\cdots,X_{\ell-1},\pm x^{-1}|x_1,\cdots,x_{2n-2},x,-x),
\quad
x=x_{n-1},
\non
\ena
and for a set $M=(m_1,\cdots,m_\ell)$, $1\leq m_1<\cdots<m_\ell\leq n$,
define
\bea
&&
M^{(a)}=(m_1',\cdots,m_{\ell-1}')=(m_1,\cdots,m_{a-1},m_{a+1},\cdots,m_\ell).
\non
\ena
Let us set
\bea
&&
\tilde{g}_M=g_M\vert_{\beta_n=\beta_{n-1}+\pi i},
\quad
D_n(X)=\prod_{j=1}^n(1-Xx_j)
\non
\\
&&
\tilde{J}_M(P)=
\prod_{b=1}^{\ell-1}\phi_{n-2}(\alpha_b)\cdot
\tilde{g}_M\cdot
\frac{P}{\prod_{b=1}^{\ell-1}D_{n-2}(X_b)(1-X_b^2x_{n-1}^2)}.
\non
\ena
Using these notations we define
\bea
&&
I^{(1)}_M(P)=
\sum_{a=1}^\ell \frac{(-1)^{\ell-1}}{\beta_{n-1}-\beta_{m_a}-\pi i}
\prod_{j=1}^{n-2}\frac{\beta_{n-1}-\beta_j}{\beta_{n-1}-\beta_j-\pi i}
\int_{C^{\ell-1}}
\prod_{b=1}^{\ell-1}d\alpha_b 
\tilde{J}_{M^{(a)}}(P)
\non
\\
&&
\times 
\frac{\prod_{b=1}^{a-1}(\alpha_b-\beta_{n-1}+2\pi i)}
{\prod_{b=a}^{\ell-1}(\beta_{n-1}-\alpha_b)},
\non
\\
&&
I^{(2)}_M(P)=\int_{C^{\ell-1}}\prod_{b=1}^{\ell-1}d\alpha_b 
\tilde{J}_{M^{(\ell)}}(P),
\non
\\
&&
I^{(3)}_M(P)=\int_{C^{\ell-1}}\prod_{b=1}^{\ell-1}d\alpha_b 
\tilde{J}_{M^{(\ell)}}(P)
\prod_{b=1}^{\ell-1}
\frac{\alpha_{b}-\beta_{n-1}}
{\alpha_{b}-\beta_{n-1}+\pi i},
\non
\\
&&
I^{(4)}_M(P)=\int_{C^{\ell-1}}\prod_{b=1}^{\ell-1}d\alpha_b 
\tilde{J}_{M^{(\ell-1)}}(P)
\frac{\beta_{n-1}-\alpha_{\ell-1}+\pi i}
{\alpha_{\ell-1}-\beta_{n-1}+\pi i},
\non
\ena
and
\bea
&&
f^{(1)}_{n-1}=\frac{-(-2\pi i)^\ell}{2}
\frac{\phi_{n-1}(\beta_{n-1}-\pi i)}{D_{n-2}(-x_{n-1}^{-1})},
\non
\\
&&
f^{(2)}_{n-1}=f^{(4)}_{n-1}=(-2\pi i)^{\ell-1}
\frac{\phi_{n-1}(\beta_{n-1})}{D_{n-2}(x_{n-1}^{-1})},
\non
\\
&&
f^{(3)}_{n-1}=(-2\pi i)^{\ell-1}
\frac{\phi_{n-1}(\beta_{n-1}+\pi i)}{D_{n-2}(-x_{n-1}^{-1})}
\prod_{j=1}^{n-2}\frac{\beta_{n-1}-\beta_j+2\pi i}{\beta_{n-1}-\beta_j}.
\non
\ena
Then

\begin{prop}\label{prop2-2}
For $P\in \wedge^{\ell}H^{(2n)}$, the residue of $I_M$ is given by
\bea
\frac{1}{2\pi i}Res_{\beta_n=\beta_{n-1}+\pi i}I_M(P)
&=&
f^{(1)}_{n-1}I^{(1)}_M(\rho_{-}(P))
\quad
\text{if $n-1,n\notin M$},
\non
\\
&=&
f^{(1)}_{n-1}I^{(1)}_M(\rho_{-}(P))
-
f^{(2)}_{n-1}I^{(2)}_M(\rho_{+}(P))
\quad
\text{if $n-1\in M$, $n\notin M$},
\non
\\
&=&
\sum_{j=1}^3f^{(j)}_{n-1}I^{(j)}_M(\rho_{(-1)^j}(P))
\quad
\text{if $n\in M$, $n-1\notin M$},
\non
\\
&=&
\sum_{j=1}^4f^{(j)}_{n-1}I^{(j)}_M(\rho_{(-1)^j}(P))
\quad
\text{if $n-1,n\in M$},
\non
\ena
where $\pm1$ are identified with $\pm$.
\end{prop}

\begin{cor}
If
\bea
&&
\rho_{+}(P)=\rho_{-}(P)=0,
\label{eq-P}
\ena
then (\ref{eq2-4}) is satisfied.
\end{cor}
\vskip2mm

\noindent
{\bf Remark}. The residue of $I_M(P)$ can vanish for some $M$ although
$\rho_+(P)$ or $\rho_-(P)$ is not zero. The $P$'s corresponding to
form factors of local operators give such examples \cite{Smir2, NT}.
We conjecture that if the residue of $I_M(P)$ vanishes for all $M$ (\ref{eq-P})
holds. In other words we conjecture that (\ref{eq2-4}) and (\ref{eq-P}) 
are equivalent. It is true for $\ell=1$.
\vskip2mm

Let us write down the equations (\ref{eq-P}) in 
component forms. 
To this end we introduce a notation.
For a polynomial $f(x_1,\cdots,x_{n})$ we define $\bar{f}$ by
\bea
&&
\bar{f}=f(x_1,\cdots,x_{n-2},x,-x), 
\quad
x=x_{n-1}.
\non
\ena
Expand $P\in \wedge^\ell H^{(2n)}$ as
\bea
&&
P=\sum\, P_{i_1 \cdots i_\ell}\, X^{i_1}\wedge\cdots\wedge X^{i_\ell},
\non
\ena
such that the coefficients $P_{i_1 \cdots i_\ell}$'s are anti-symmetric
with respect to indices:
\bea
&&
P_{i_{\sigma_1} \cdots i_{\sigma_\ell}}=
\sgn\, \sigma\cdot P_{i_1 \cdots i_\ell},
\quad
\sigma\in S_\ell.
\non
\ena
Then (\ref{eq-P}) is equivalent to 
\bea
&&
\sum_{i_\ell:\, even}\overline{ P_{i_1 \cdots i_\ell}}x^{-i_\ell}
=
\sum_{i_\ell:\, odd}\overline{ P_{i_1 \cdots i_\ell}}x^{-i_\ell}
=0.
\label{eq-P-comp}
\ena

Let $U_{n,\ell}$ be the space of solutions of (\ref{eq-P}):
\bea
&&
U_{n,\ell}=\{\, P\in \wedge^\ell H^{(n)}\vert \, 
\rho_{+}(P)=\rho_{-}(P)=0\,\}.
\non
\ena
It is obviously an $R_n$-module. In the subsequent sections we shall
study the structure of $U_{n,\ell}$.

\section{Solutions for $\ell=1$}
In this section we study $U_{2n,1}$.
The equations (\ref{eq-P-comp}) for $P=\sum_{i=0}^{2n-1}P_i(x)X^i$ are
\bea
&&
\sum_{i=0}^{n-1}\overline{P_{2i}(x)}x^{-2i}=
\sum_{i=1}^{n}\overline{P_{2i-1}(x)}x^{-(2i-1)}=0.
\label{eqevenl=1}
\ena
Thus the sets of polynomials $\{P_{2i}\}$ and $\{P_{2i-1}\}$ satisfy 
the same equation independently.
The solutions of the equation (\ref{eqevenl=1}) can be found easily 
as we shall see below.

Let $e^{(n)}_k$ be the elementary symmetric polynomial defined by
\bea
&&
\prod_{j=1}^n(1+x_jt)=\sum_{k=0}^ne^{(n)}_kt^k.
\non
\ena
We set $e^{(n)}_k=0$ for $k>n$, $k<0$ or $n<0$ and $e^{(0)}_0=1$.
\vskip2mm

\noindent
Example 3.1. The polynomial
\bea
&&
P=e^{(2n)}_1X+ e^{(2n)}_3X^3+\cdots+ e^{(2n)}_{2n-1}X^{2n-1}
\non
\ena
is a solution of (\ref{eqevenl=1}). To verify this we use
\bea
&&
\overline{e^{(n)}_k}=e^{(n-2)}_k-x^2e^{(n-2)}_{k-2},
\label{recursionE}
\ena
which is valid for all integers $k$. Then
\bea
\sum_{s=1}^n\overline{e^{(2n)}_{2s-1}}x^{-(2s-1)}
&=&
\sum_{s=1}^n(e^{(2n-2)}_{2s-1}-x^2e^{(2n-2)}_{2s-3})x^{-(2s-1)}
\non
\\
&=&\sum_{s=1}^ne^{(2n-2)}_{2s-1}x^{-(2s-1)}-
\sum_{s=1}^ne^{(2n-2)}_{2s-1} x^{-(2s-1)}=0.
\non
\ena 
\vskip2mm

Beginning from the solution in Example 3.1 we can construct a set of solutions of (\ref{eqevenl=1}).
Let $P^{(2n)}_{r,s}$, $1\leq r\leq n$, $s\in \mathbb{Z}$ be defined
by the following recursion relations:
\bea
P^{(2n)}_{1,s}&=&e^{(2n)}_{2s-1},
\non
\\
P^{(2n)}_{r,s}&=&P^{(2n)}_{r-1,s+1}-e^{(2n)}_{2s}P^{(2n)}_{r-1,1}
\quad \hbox{for $r\geq 2$}.
\label{def-P}
\ena
We easily have
\bea
&&
P^{(2n)}_{r,s}=0,
\quad
s\leq 0.
\label{s-negative}
\ena

\begin{prop}\label{soll=1} We have 
\bea
&&
\sum_{s=1}^n\overline{P^{(2n)}_{r,s}}x^{-2(s-1)}=0.
\non
\ena 
\end{prop}
\vskip2mm

\noindent
{\it Proof.} 
Let us prove this equation by the induction on $r$.
The case $r=1$ is verified in Example 3.1.

\begin{lemma}\label{lem0}
$P^{(2n)}_{r,s}=0$, $s\geq n+1$.
\end{lemma}

This lemma is easily proved by the induction on $r$.

\begin{lemma}\label{lem1}
$\overline{P^{(2n)}_{r,s}}=P^{(2n-2)}_{r,s}-x^2P^{(2n-2)}_{r,s-1}$.
\end{lemma}
\vskip2mm

\noindent
{\it Proof.} Let us prove the lemma by the induction on $r$.
If $r=1$, the assertion is (\ref{recursionE}) and there is nothing to prove.
For $r\geq 2$, using (\ref{def-P}) and (\ref{recursionE}) we have
\bea
\overline{P^{(2n)}_{r,s}}
&=&
\overline{P^{(2n)}_{r-1,s+1}}-
\overline{e^{(2n)}_{2s}}\overline{P^{(2n)}_{r-1,1}}
\non
\\
&=&
P^{(2n-2)}_{r-1,s+1}-x^2P^{(2n-2)}_{r-1,s}
-(e^{(2n-2)}_{2s}-x^2e^{(2n-2)}_{2s-2})P^{(2n-2)}_{r-1,1}
\non
\\
&=&
P^{(2n-2)}_{r-1,s+1}-e^{(2n-2)}_{2s}P^{(2n-2)}_{r-1,1}
-x^2(P^{(2n-2)}_{r-1,s}-e^{(2n-2)}_{2s-2}P^{(2n-2)}_{r-1,1})
\non
\\
&=&
P^{(2n-2)}_{r,s}-x^2P^{(2n)}_{r,s-1}.
\non
\ena
Q.E.D.
\vskip2mm

Using Lemma \ref{lem0}, \ref{lem1} and (\ref{s-negative}) we have, 
for $r\geq 2$,
\bea
\sum_{s=1}^n\overline{P^{(2n)}_{r,s}}x^{-2(s-1)}
&=&
\sum_{s=1}^n(P^{(2n-2)}_{r,s}-x^2P^{(2n-2)}_{r,s-1})x^{-2(s-1)}=0.
\non
\ena
Thus Proposition \ref{soll=1} is proved.
Q.E.D.
\vskip2mm

Define the elements $v^{(2n)}_r$, $w^{(2n)}_r$ of $H^{(2n)}$ by
\bea
&&
v^{(2n)}_r=\sum_{s=1}^nP^{(2n)}_{r,s}X^{2(s-1)},
\qquad
w^{(2n)}_r=Xv^{(2n)}_r.
\label{def-vw}
\ena

\begin{cor}
The elements $v^{(2n)}_r$, $w^{(2n)}_r$ are solutions of (\ref{eqevenl=1}).
\end{cor}

\begin{theorem}\label{un1}
The space $U_{2n,1}$ is a free $R_{2n}$ module of rank $2n$ with 
$\{v^{(2n)}_r, w^{(2n)}_r\}$ as a basis:
$$
U_{2n,1}=\oplus_{r=1}^n R_{2n} v^{(2n)}_r\oplus_{r=1}^n R_{2n}w^{(2n)}_r.
$$
\end{theorem}
\vskip2mm

For the proof let us define the $n$ by $n$ matrix $X^{(2n)}=(X^{(2n)}_{ij})$ by
\bea
&&
X^{(2n)}_{ij}=P^{(2n)}_{i,j}, \quad 1\leq i,j\leq n.
\non
\ena

\begin{prop}\label{detx2n}
Let $\Delta^{+}_n=\prod_{i<j}^n(x_i+x_j)$. Then
$$
\det X^{(2n)}=\Delta^{+}_{2n}.
$$
\end{prop}
\vskip2mm

We prepare several lemmas in order to prove Proposition \ref{detx2n}.

\begin{lemma}\label{lem2}
As a polynomial of $x_1$, ..., $x_{2n}$, 
$\det X^{(2n)}$ is homogeneous of degree 
$\binom{2n}{2}$.
\end{lemma}
\vskip2mm

\noindent
\pf It can be easily checked that $X^{(2n)}_{ij}=P^{(2n)}_{i,j}$ 
is homogeneous of degree $2i+2j-3$. Thus for any permutation 
$(j_1,...,j_n)$ of $(1,...,n)$ we have
\bea
&&
\deg(X^{(2n)}_{1j_1}\cdots X^{(2n)}_{nj_n})=2\sum_{i=1}^ni+2\sum_{j=1}^nj-3n=
n(2n-1).
\non
\ena
Q.E.D.

\begin{lemma}\label{lem3}
$\det X^{(2n)}$ is divisible by $\Delta^{+}_{2n}$.
\end{lemma}
\vskip2mm

\noindent
\pf Since $\det X^{(2n)}$ is a symmetric polynomial, it is sufficient to prove
\bea
&&
\overline{\det X^{(2n)}}=0.
\label{detx2n0}
\ena
In $\det X^{(2n)}$ we successively add the $j$-th column multiplied by 
$x_{2n-1}^{-2(j-1)}$ to the first column. Then $i$-th row of the first column
becomes 
\bea
&&
\sum_{j=1}^n P^{(2n)}_{i,j}x^{-2(j-1)}.
\non
\ena
Thus  (\ref{detx2n0}) follows from Proposition \ref{soll=1}.
Q.E.D.

\begin{lemma}\label{lem4}
The expression of $\det X^{(2n)}$ as a polynomial of $e_1$,...,$e_{2n}$
is
\bea
&&
\det X^{(2n)}=(-1)^{\frac{1}{2}n(n-1)}e_1^ne_{2n}^{n-1}+\cdots,
\non
\ena
where $e_k=e^{(2n)}_k$ and $\cdots$ part does not contain $e_1^ne_{2n}^{n-1}$.
\end{lemma}
\vskip2mm

\noindent
\pf It can be easily proved that $P^{(2n)}_{r,s}$ is at most of degree $1$
in $e_{2n}$, is of degree $0$ for $s<n-r+2$ and
\bea
&&
P^{(2n)}_{r,n-r+2}=-e_1e_{2n}+(\hbox{terms of degree $0$ in $e_{2n}$}).
\non
\ena
Thus
$$
\lim_{e_{2n}\rightarrow \infty}
\frac{1}{e_{2n}^{n-1}}\det X^{(2n)}
=
\left|
\begin{array}{rrrrr}
e_1& e_3&\cdots&e_{2n-3}&e_{2n-1}\\
0&0&\cdots&0&-e_1\\
0&0&\cdots&-e_1&\ast\\
\vdots&\vdots&{}&{}&\vdots\\
0&-e_1&\cdots&\cdots&\ast\\
\end{array}
\right|
=
(-1)^{\frac{1}{2}n(n-1)}e_1^n,
$$
which proves the lemma.
Q.E.D.

\begin{lemma}\label{lem5}
The expression of $\Delta^{+}_{2n}$ as a polynomial of $e_1$,...,$e_{2n}$
is
\bea
&&
\Delta^{+}_{2n}=(-1)^{\frac{1}{2}n(n-1)}e_1^ne_{2n}^{n-1}+\cdots.
\non
\ena
\end{lemma}
\vskip2mm

\noindent
\pf The statement follows easily from (see \cite{Mac} for example)
\bea
&&
\Delta^{+}_{2n}=\det(e^{(2n)}_{2n-2i+j})_{1\leq i,j\leq 2n-1}.
\non
\ena
Q.E.D.
\vskip3mm

Now Proposition \ref{detx2n} is obvious from Lemma \ref{lem3}, \ref{lem4} and
\ref{lem5}.
Q.E.D.
\vskip4mm

\noindent
{\it Proof of Theorem \ref{un1}.}
\par
\noindent
By Proposition \ref{detx2n}, $v_i^{(2n)}$, $w_j^{(2n)}$, $i,j=1,...,n$ are 
linearly independent over $R_{2n}$. Thus it is sufficient to prove that 
any solution
$$
P=\sum_{s=1}^n P_s X^{2(s-1)}
$$ 
of (\ref{eqevenl=1}) is written as a linear combination of 
$v_i^{(2n)}$, $1\leq i\leq n$ with the coefficients in $R_{2n}$.
Let us consider the following linear equation for $Q_i$, $i=1,...,n$:
\bea
&&
P=\sum_{i=1}^n Q_i v_i^{(2n)},
\non
\ena
which, in the matrix form, is
\bea
&&
(P_1,\cdots,P_n)=(Q_1,\cdots,Q_n)X^{(2n)}.
\non
\ena
By Cramer's formula, using Proposition \ref{detx2n}, this equation is 
solved as
\bea
&&
Q_i=\frac{1}{\Delta^{+}_{2n}}\det X^{(2n)}_i,
\non
\ena
where the matrix $X^{(2n)}_i$ is defined by replacing $i$-th row in $X^{(2n)}$
by $(P_1,\cdots,P_n)$.
In a similar way to Lemma \ref{lem3} one can prove that $X^{(2n)}_i$
is divisible by $\Delta^{+}_{2n}$. Thus $Q_i\in R_{2n}$ and the theorem
is proved.
Q.E.D.

Notice that if $\psi_P=0$ the residue of $\psi_P$ is zero.
Thus $\Xi_1^{(2n)}$ should be in $U_{2n,\ell}$.
In fact we have the relation
\bea
&&
\Xi_1^{(2n)}=w_1^{(2n)},
\label{w1=Xi1}
\ena
which is obvious from their definitions.

\section{Solutions for $\ell=2$}
In this section we solve the equations (\ref{eq-P}) in the case $\ell=2$.
In general, for solutions $u_j$, $1\leq j\leq m$, of (\ref{eq-P}) with 
$\ell=\ell_j$, it is obvious that their exterior product
\bea
&&
u_1\wedge \cdots\wedge u_m=
\sum_\sigma \sgn\, \sigma\, 
u_1(X_{\sigma_1},\cdots,X_{\sigma_{\ell_1}})\cdots 
u_m(X_{\sigma_{\ell_1+\cdots+\ell_{m-1}+1}},\cdots,X_{\sigma_\ell}),
\non
\ena
is a solution of (\ref{eq-P}) with $\ell=\ell_1+\cdots+\ell_m$.
Thus, from the results of the previous section, 
\bea
&&
v^{(2n)}_{i_1}\wedge v^{(2n)}_{i_2},
\quad
v^{(2n)}_{i}\wedge w^{(2n)}_{j},
\quad
w^{(2n)}_{j_1}\wedge w^{(2n)}_{j_2},
\non
\ena
are solutions of (\ref{eq-P}) with $\ell=2$ for any $i_1,i_2$, $i,j$, 
$j_1,j_2$.
In fact there exist solutions which can not be expressed as an 
$R_{2n}$-linear combination of those solutions. 
Let us give a construction of them.

Set 
\bea
&&
v_0^{(2n)}=\sum_{j=0}^{n}e^{(2n)}_{2j}X^{2j}.
\non
\ena
Although $v_0^{(2n)}$ are not in $H^{(2n)}$ it satisfies
\bea
&&
\rho_\pm(v_0^{(2n)})=0.
\non
\ena
Using $v_0^{(2n)}$ we define, for $1\leq j\leq n$,
\bea
2\xi_j^{(2n)}(X_1,X_2)&=&
\frac{X_1-X_2}{X_1+X_2}
\Big(
v_0^{(2n)}(X_1)w^{(2n)}_j(X_2)+v_0^{(2n)}(X_2)w^{(2n)}_j(X_1)
\Big)
\non
\\
&&
-v_0^{(2n)}(X_1)w^{(2n)}_j(X_2)+v_0^{(2n)}(X_2)w^{(2n)}_j(X_1).
\label{def-xi}
\ena
Since 
$$
v_0^{(2n)}(-X)=v_0^{(2n)}(X),
\quad
w^{(2n)}_j(-X)=-w^{(2n)}_j(X),
$$
$\xi_j$ is an anti-symmetric polynomial of $X_1$ and $X_2$ with the
coefficients in $R_{2n}$.
The definition of $\xi_j$ is motivated by the relation
\bea
&&
\Xi_2^{(2n)}=\xi_1^{(2n)},
\label{Xi2=xi1}
\ena
which can be easily checked. 

\begin{prop}
$\xi_j^{(2n)}\in U_{2n,2}, \quad 1\leq j\leq n$.
\end{prop}
\vskip2mm

\noindent
\pf Applying $\rho_\pm$ to (\ref{def-xi}) we see that
\bea
\rho_\pm(2\xi_j^{(2n)})&=&
\frac{X_1\mp x^{-1}}{X_1\pm x^{-1}}
\Big(
\overline{v_0^{(2n)}(X_1)}\rho_\pm(w^{(2n)}_j)+
\rho_\pm(v_0^{(2n)})\overline{w^{(2n)}_j(X_1)}
\Big)
\non
\\
&&
-
\overline{v_0^{(2n)}(X_1)}\rho_\pm(w^{(2n)}_j)+
\rho_\pm(v_0^{(2n)})\overline{w^{(2n)}_j(X_1)}
\non
\\
&=&0,
\non
\ena
where $\overline{v_0^{(2n)}(X_1)}$ means, for example, that
we take bar of each coefficient of the power of $X_1$.

It remains to prove $\deg_{X_i} \xi_j^{(2n)}\leq 2n-1$, $i=1,2$.
Obviously $\deg_{X_i} \xi_j^{(2n)}\leq 2n$, $i=1,2$. 
Let us calculate the coefficient
of $X_1^{2n}$ as a polynomial of $X_1$ as
\bea
&&
\lim_{X_1\rightarrow \infty}\frac{2\xi_j^{(2n)}}{X_1^{2n}}
\non
\\
&=&
\lim_{X_1\rightarrow \infty}
\Big[
\frac{ 1-\frac{X_2}{X_1} }{ 1+\frac{X_2}{X_1} }
\Big(
\frac{ v_0^{(2n)}(X_1) }{ X_1^{2n} } w^{(2n)}_j(X_2)
+
v_0^{(2n)}(X_2) \frac{ w^{(2n)}_j(X_1) }{ X_1^{2n} }
\Big)
\non
\\
&&
-\frac{ v_0^{(2n)}(X_1) }{ X_1^{2n} } w^{(2n)}_j(X_2)
+v_0^{(2n)}(X_2) \frac{ w^{(2n)}_j(X_1) }{ X_1^{2n} }
\Big]
\non
\\
&=&
e_{2n}^{(2n)}w^{(2n)}_j(X_2)-e_{2n}^{(2n)}w^{(2n)}_j(X_2)=0.
\non
\ena
Thus $\deg_{X_1} \xi_j^{(2n)}\leq 2n-1$. By the anti-symmetry of $\xi_j^{(2n)}$
we get $\deg_{X_2} \xi_j^{(2n)}\leq 2n-1$.
Q.E.D.

\begin{theorem}\label{un2}
The module $U_{2n,2}$ is a free $R_{2n}$ module of rank 
$\textstyle{\binom{2n}{2}}$
with the following elements as a basis
\bea
&&
v_{i_1}^{(2n)}\wedge v_{i_2}^{(2n)}\,\,(1\leq i_1<i_2\leq n),
\quad
w_{j_1}^{(2n)}\wedge w_{j_2}^{(2n)}\,\,(1\leq j_1<j_2\leq n),
\non
\\
&&
v_i^{(2n)}\wedge w_j^{(2n)}\,\,(1\leq i\leq n, 1\leq j\leq n-1),
\quad
\xi_k^{(2n)}\,\, (1\leq k\leq n).
\label{basis-un2}
\ena
\end{theorem}
\vskip2mm

The linear independence of (\ref{basis-un2}) is a special case of
Theorem \ref{conj-indep}.
Let us prove the theorem assuming the linear independence.

First of all one can easily check that the number of elements 
(\ref{basis-un2}) is $\binom{2n}{2}$.
We number the elements (\ref{basis-un2}) from $1$ to $\binom{2n}{2}$ and
name the $r$-th element by $Q_r$. 
Let 
\bea
&&
Q_r=\sum_{i<j}Q_{r;ij}X^i\wedge X^j,
\quad
Q_{r;ij}\in R_{2n}.
\non
\ena
We set $Q_{r;ij}=-Q_{r;ji}$ for $i\geq j$.
We introduce the lexicographical order from the left on the set
$$
J_2=\{\,(ij)\,\vert\,\, 0\leq i<j\leq 2n-1\,\}.
$$

Define the $\binom{2n}{2}$ by $\binom{2n}{2}$ 
matrix $X^{(2n)}=(X^{(2n)}_{r;ij})$,
$1\leq r\leq \binom{2n}{2}$, $(ij)\in J_2$ by
\bea
&&
X^{(2n)}_{r;ij}=Q_{r;ij}.
\non
\ena

\begin{lemma}\label{lem21}
As a polynomial of $x_1$, ..., $x_{2n}$, $\det X^{(2n)}$ is homogeneous 
of degree $n(2n-1)(4n-3)$.
\end{lemma}
\vskip2mm

This is a special case of Proposition \ref{degX}.

\begin{lemma}\label{lem22}
$\det X^{(2n)}$ is divisible by $(\Delta^{+}_{2n})^{4n-3}$.
\end{lemma}
\vskip1mm

This lemma is also a special case of Proposition \ref{divisible} below.
Since the proof of the lemma generalizes the case $\ell=1$ and
illustrates the proof in the general case, we shall give a proof of it here.
\vskip2mm

\noindent
\pf In $\det X^{(2n)}$, to each $(0,j)$-th column with $j=1,...,2n-1$, 
we successively add the $(2i,j)$-th column times $x_{2n-1}^{-2i}$ for 
$2\leq 2i<j$ and the $(j,2i)$-th column times $-x_{2n-1}^{-2i}$ for $j<2i$.
Then, to each $(1,j)$-th column with $j=2,...,2n-1$, we add 
the $(2i-1,j)$-th column times $x_{2n-1}^{-(2i-1)}$ for $3\leq 2i-1<j$
and the $(j,2i-1)$-th column times $-x_{2n-1}^{-(2i-1)}$ for 
$j<2i-1$ successively.
Then, in the resulting determinant, the $(r;0,j)$-th component $Q'_{r;0,j}$
and the $(r;1,j')$-th component $Q'_{r;1,j'}$ becomes
\bea
&&
Q'_{r;0,j}=\sum_{i=0}^{n-1}Q_{r;2i,j}x^{-2i}_{2n-1},
\quad
Q'_{r;1,j'}=\sum_{i=1}^{n}Q_{r;2i-1,j'}x^{-(2i-1)}_{2n-1},
\non
\ena
respectively. Since $Q_r$ is a solution of (\ref{eq-P-comp}),
\bea
&&
\overline{Q'_{r;0,j}}=\overline{Q'_{r;1,j'}}=0.
\non
\ena
This means that $Q'_{r;0,j}$ and $Q'_{r;1,j'}$ are divisible by 
$x_{2n-1}+x_{2n}$. Thus $\det X^{(2n)}$ is divisible by 
$(x_{2n-1}+x_{2n})^{4n-3}$. Because $\det X^{(2n)}$ is a symmetric polynomial
of $x_1,...,x_{2n}$, it is divisible by $(\Delta^{+}_{2n})^{4n-3}$.
Q.E.D.
\vskip2mm

The following lemma is obvious from the linear independence of 
(\ref{basis-un2}) and Lemma \ref{lem21}, \ref{lem22}.

\begin{lemma}\label{lem23}
There is a non-zero constant $c$ such that
\bea
&&
\det X^{(2n)}=c(\Delta^{+}_{2n})^{4n-3}.
\non
\ena
\end{lemma}
\vskip2mm

The proof of Theorem \ref{un2} is similar to that of Theorem \ref{un1}
using Lemma \ref{lem2} and \ref{lem3}.

\section{$q$-Tetranomial identity}\label{qtetra-sec}
From now on we use the multi-index notations like
\bea
v_I=v_{i_1}\wedge\cdots\wedge v_{i_l},
\quad
\text{for $I=(i_1,\cdots,i_l)$}.
\non
\ena

Consider the following elements of $U_{2n,\ell}$:
\bea
&&
v^{(2n)}_I\wedge w^{(2n)}_J\wedge \xi^{(2n)}_K,
\label{basis-unl}
\\
&&
I=(i_1,\cdots,i_{\ell_1}),
\quad
1\leq i_1<\cdots<i_{\ell_1}\leq n,
\label{indexI}
\\
&&
J=(j_1,\cdots,j_{\ell_2}),
\quad
1\leq j_1<\cdots<j_{\ell_2}\leq n-\ell_1-\ell_3,
\label{indexJ}
\\
&&
K=(k_1,\cdots,k_{\ell_3}),
\quad
1\leq k_1\leq \cdots \leq k_{\ell_3}\leq n-\ell_1-\ell_3+1,
\label{indexK}
\\
&&
\ell_1+\ell_2+2\ell_3=\ell.
\label{sum=l}
\ena
If $\ell=1,2$, these are the basis of $U_{2n,1}$ and $U_{2n,2}$
given in Theorem \ref{un1} and \ref{un2} respectively.

Let us count the number $N_{2n,\ell}$ of such elements.
Obviously it is given by
\bea
&&
N_{2n,\ell}=
\sum_{\ell_1+\ell_2+2\ell_3=\ell}
\binom{n}{\ell_1}\binom{n-\ell_1}{\ell_3}\binom{n-\ell_1-\ell_3}{\ell_2}.
\non
\ena
In the right hand side
\bea
&&
\binom{n}{\ell_1}\binom{n-\ell_1}{\ell_3}\binom{n-\ell_1-\ell_3}{\ell_2}=
\frac{n!}{\ell_1!\ell_2!\ell_3!(n-\ell_1-\ell_2-\ell_3)!},
\non
\ena
which is the tetranomial coefficient denoted by
\bea
&&
\tetc{n}{\ell_1}{\ell_2}{\ell_3}.
\label{def-tetra}
\ena
Therefore we have
\bea
&&
(1+x+x+x^2)^n=\sum_{\ell=0}^{2n}N_{2n,\ell}x^\ell.
\non
\ena
The left hand side is
\bea
&&
(1+x)^{2n}=\sum_{\ell=0}^{2n}\binom{2n}{\ell}x^\ell.
\non
\ena
Thus we have
\bea
&&
N_{2n,\ell}=\binom{2n}{\ell}=\rk\, \wedge^\ell H^{(2n)},
\label{dim-correct}
\ena
and the identity
\bea
&&
\binom{2n}{\ell}=
\sum_{\ell_1+\ell_2+2\ell_3=\ell}
\binom{n}{\ell_1}\binom{n-\ell_1}{\ell_2}\binom{n-\ell_1-\ell_2}{\ell_3},
\label{tetra-id}
\ena
here we used the symmetry of (\ref{def-tetra}) with respect to
the permutations of $\ell_1$, $\ell_2$, $\ell_3$.

By the definition of $\deg_1$ we easily have
\bea
&&
\deg_1\, v^{(2n)}_i=2i-1,
\quad
\deg_1\, w^{(2n)}_j=2j-2,
\quad
\deg_1\, \xi^{(2n)}_k=2k-2.
\non
\ena
Thus 
\bea
&&
\deg_1\,(\,v^{(2n)}_I\wedge w^{(2n)}_J\wedge \xi^{(2n)}_K\,)
=
\sum_{r=1}^{\ell_1}(2i_r-1)
+\sum_{r=1}^{\ell_2}(2j_r-2)
+\sum_{r=1}^{\ell_3}(2k_r-2),
\label{deg1}
\ena
if the element is non-zero. We denote the right hand side of (\ref{deg1})
by $d_{I,J,K}$.
Let us consider the generating function of $d_{I,J,K}$:
\bea
&&
\ch_{2n,\ell}:=\sum_{I,J,K}q^{d_{I,J,K}},
\label{ch2nl-1}
\ena
where $I,J,K$ run the index sets satisfying (\ref{indexI}) - (\ref{sum=l}).
It is easy to evaluate it explicitly. To write it we use the 
$q$-integer notations:
\bea
&&
[n]_p=1-p^n,
\quad
[n]_p!=[1]_p[2]_p\cdots[n]_p,
\quad
\pbc{n}{r}=\frac{[n]_p!}{[r]_p![n-r]_p!}.
\non
\ena
Then we easily obtain
\bea
&&
\ch_{2n,\ell}=
\sum_{\ell_1+\ell_2+2\ell_3=\ell}
q^{\ell_1^2+\ell_2(\ell_2-1)}
\qqbc{n}{\ell_1}\qqbc{n-\ell_1}{\ell_3}\qqbc{n-\ell_1-\ell_3}{\ell_2}.
\label{ch2nl-2}
\ena
In the right hand side 
\bea
\qqbc{n}{\ell_1}\qqbc{n-\ell_1}{\ell_2}\qqbc{n-\ell_1-\ell_2}{\ell_3}
=
\frac{[n]_{q^2}!}{[\ell_1]_{q^2}![\ell_2]_{q^2}![\ell_3]_{q^2}!
[n-\ell_1-\ell_2-\ell_3]_{q^2}!}
\non
\ena
is the $q$-tetranomial coefficient written as
\bea
&&
\qqtetc{n}{\ell_1}{\ell_2}{\ell_3}.
\label{qtetra}
\ena
If some $\ell_i\leq -1$, we set $\text{(\ref{qtetra})}=0$.
Notice that (\ref{qtetra}) is symmetric with 
respect to the permutations of $\ell_1$, $\ell_2$, $\ell_3$.
Now we prove the $q$-analogue of (\ref{tetra-id}).

\begin{prop}\label{q-tetra-id}
\bea
&&
\qbc{2n}{\ell}
=
\sum_{\ell_1+\ell_2+2\ell_3=\ell}
q^{\ell_1^2+\ell_2(\ell_2-1)}
\qqbc{n}{\ell_1}\qqbc{n-\ell_1}{\ell_2}\qqbc{n-\ell_1-\ell_2}{\ell_3}.
\label{qtetra-id}
\ena
\end{prop}
\vskip2mm
\noindent
\pf Let us denote the left and the right hand sides of (\ref{qtetra-id}) by
$a_{n,\ell}$ and $b_{n,\ell}$ respectively.
By direct calculations one can easily prove (\ref{qtetra-id}) for $\ell=0,1,2$.
It is obvious that $a_{1,\ell}=b_{1,\ell}=0$ for $\ell\geq 3$.
Thus (\ref{qtetra-id}) holds for $n=1$. Let us prove the identity
by the induction on $n$. To this end we shall show that 
$a_{n,\ell}$ and $b_{n,\ell}$ satisfy the same recursion relation.

The following recursion relations hold for $q$-binomial coefficients:
\bea
\pbc{m}{r}&=&\pbc{m-1}{r}+p^{m-r}\pbc{m-1}{r-1}
\label{binom-rec1}
\\
&=& p^r\pbc{m-1}{r}+\pbc{m-1}{r-1}.
\label{binom-rec2}
\ena
Using (\ref{binom-rec1}) and (\ref{binom-rec2}) we have
\bea
&&
a_{n,\ell}=a_{n-1,\ell}+(q^{2n-\ell-1}+q^{2n-1})a_{n-1,\ell-1}
+q^{2n-\ell}a_{n-1,\ell-2}.
\label{a-rec}
\ena

Let us prove the same recursion relation for $b_{n,\ell}$.
Using (\ref{binom-rec1}) for the first $q$-binomial coefficient in the product
in $b_{n,\ell}$ we have
\bea
&&
b_{n,\ell}
=
\sum_{\ell_1+\ell_2+2\ell_3=\ell}
q^{\ell_1^2+\ell_2(\ell_2-1)}
\qqbc{n-1}{\ell_1}\qqbc{n-\ell_1}{\ell_2}\qqbc{n-\ell_1-\ell_2}{\ell_3}
\label{b1}
\\
&&+
q^{2n-1}
\sum_{\ell_1+\ell_2+2\ell_3=\ell}
q^{(\ell_1-1)^2+\ell_2(\ell_2-1)}
\qqbc{n-1}{\ell_1-1}\qqbc{n-\ell_1}{\ell_2}\qqbc{n-\ell_1-\ell_2}{\ell_3}.
\label{b2}
\ena
It is obvious that 
\bea
&&
\text{(\ref{b2})}=q^{2n-1}b_{n-1,\ell-1}.
\non
\ena
Applying (\ref{binom-rec1}) to the second $q$-binomial coefficient 
in (\ref{b1}) we have
\bea
&&
\text{(\ref{b1})}=
\sum_{\ell_1+\ell_2+2\ell_3=\ell}
q^{\ell_1^2+\ell_2(\ell_2-1)}
\qqbc{n-1}{\ell_1}\qqbc{n-\ell_1-1}{\ell_2}\qqbc{n-\ell_1-\ell_2}{\ell_3}
\label{b3}
\\
&&
+\sum_{\ell_1+\ell_2+2\ell_3=\ell-1}\!\!\!\!\!\!\!\!
q^{\ell_1^2+\ell_2(\ell_2+1)+2(n-1-\ell_1-\ell_2)}\!\!
\qqtetc{n-1}{\ell_1}{\ell_2}{\ell_3}.
\label{b4}
\ena
We again use (\ref{binom-rec1}) to the third product of 
$q$-binomial coefficients in (\ref{b3}) we get
\bea
&&
\text{(\ref{b3})}=b_{n-1,\ell}
\non
\\
&&
+\!\!\!\!\!\!\!\!\!
\sum_{\ell_1+\ell_2+2\ell_3=\ell-2}\!\!\!\!\!\!\!
q^{\ell_1^2+\ell_2(\ell_2-1)+2n-\ell-(\ell_1+\ell_2)}\!\!
\qqtetc{n-1}{\ell_1}{\ell_2}{\ell_3}.
\label{b5}
\ena
Therefore we have
\bea
&&
b_{n,\ell}=b_{n-1,\ell}+q^{2n-1}b_{n-1,\ell-1}+
\text{(\ref{b4})}+\text{(\ref{b5})}.
\label{brec-1}
\ena

Using the symmetry of $\ell_1$ and $\ell_2$ of the $q$-tetranomial coefficient
we see that
\bea
&&
\text{(\ref{b4})}=
\sum_{\ell_1+\ell_2+2\ell_3=\ell-1}
q^{\ell_1^2+\ell_2(\ell_2-1)+2\ell_3+2n-1-\ell}
\qqtetc{n-1}{\ell_1}{\ell_2}{\ell_3}.
\label{b6}
\ena
There is an obvious identity of the form
\bea
&&
[\ell_3]_{q^2}
\qqtetc{n}{\ell_1}{\ell_2}{\ell_3}
=
[\ell_1+1]_{q^2}
\qqtetc{n}{\ell_1+1}{\ell_2}{\ell_3-1},
\label{qtetra-id2}
\ena
which is also valid for $\ell_3=0$.
We rewrite (\ref{qtetra-id2}) as
\bea
q^{2\ell_3}\qqtetc{n}{\ell_1}{\ell_2}{\ell_3}
&=&
\qqtetc{n}{\ell_1}{\ell_2}{\ell_3}
-
\qqtetc{n}{\ell_1+1}{\ell_2}{\ell_3-1}
\non
\\
&&
+
q^{2(\ell_1+1)}
\qqtetc{n}{\ell_1+1}{\ell_2}{\ell_3-1}.
\label{qtetra-id3}
\ena
We substitute (\ref{qtetra-id3}) into (\ref{b6}) and get
\bea
&&
\text{(\ref{b4})}=
q^{2n-1-\ell} b_{n-1,\ell-1}
\non
\\
&&
-
\sum_{\ell_1+\ell_2+2\ell_3=\ell-1}
q^{\ell_1^2+\ell_2(\ell_2-1)+2n-1-\ell}
\qqtetc{n-1}{\ell_1+1}{\ell_2}{\ell_3-1}
\label{b7}
\\
&&
+
\sum_{\ell_1+\ell_2+2\ell_3=\ell-1}
q^{\ell_1^2+\ell_2(\ell_2-1)+2n-1-\ell+2(\ell_1+1)}
\qqtetc{n-1}{\ell_1+1}{\ell_2}{\ell_3-1}.
\label{b8}
\ena
If $\ell_1=-1$, the sum of the summand of (\ref{b7}) and (\ref{b8})
becomes zero. Therefore the sum in (\ref{b7}) and (\ref{b8}) can be taken
for $\ell_1\geq -1$, $\ell_2\geq 0$, $\ell_3\geq 1$.
Then shifting the summation variables $\ell_1$ and $\ell_3$ we have
\bea
&&
\text{(\ref{b7})}+ \text{(\ref{b8})}
=-\text{(\ref{b5})}+q^{2n-\ell}b_{n-1,\ell-2}.
\non
\ena
We substitute this into (\ref{brec-1}) and get the recursion relation 
(\ref{a-rec}) for $b_{n,\ell}$.
Thus by the induction on $n$ Proposition \ref{q-tetra-id} is proved.
Q.E.D.

\section{Solutions for general $\ell$}
We shall show that the solutions of (\ref{eq-P}) for $\ell\geq 3$ 
are generated by the solutions for $\ell=1,2$. More precisely we prove

\begin{theorem}\label{conj-basis}
The module $U_{2n,\ell}$ is a free $R_{2n}$-module of rank $\binom{2n}{\ell}$
with the set of elements $\{v_I^{(2n)}\wedge w_J^{(2n)}\wedge \xi_K^{(2n)}\}$
as a basis:
\bea
&&
U_{2n,\ell}=
\oplus_{I,J,K}R_{2n} v_I^{(2n)}\wedge w_J^{(2n)}\wedge \xi_K^{(2n)},
\non
\ena
where $I,J,K$ runs over all index sets satisfying (\ref{indexI})-(\ref{sum=l}).
\end{theorem}
\vskip2mm

The strategy of the proof is similar to the case $\ell=1,2$.
The most complex part of the proof is in the following theorem.

\begin{theorem}\label{conj-indep}
The elements $\{v_I^{(2n)}\wedge w_J^{(2n)}\wedge \xi_K^{(2n)}\}$,
where $I,J,K$ runs over all index sets satisfying 
(\ref{indexI})-(\ref{sum=l}), are linearly independent over $R_{2n}$.
\end{theorem}

The proof of this theorem is given in section \ref{pf-indep}.
In this section we prove Theorem \ref{conj-basis} assuming 
Theorem \ref{conj-indep}.

We number the elements $\{v_I^{(2n)}\wedge w_J^{(2n)}\wedge \xi_K^{(2n)}\}$ 
from $1$ to $\binom{2n}{\ell}$ and denote the $r$-th element by $Q_r$.
Expand $Q_r$ as
\bea
&&
Q_r=
\sum_{0\leq i_1<\cdots< i_\ell\leq 2n-1}
Q_{r;i_1\cdots i_\ell}X^{i_1}\wedge\cdots\wedge X^{i_\ell}.
\label{qr}
\ena
We extend the index $(i_1,\cdots,i_\ell)$ of $Q_{r;i_1\cdots i_\ell}$
to $0\leq i_1,\cdots,i_\ell\leq 2n-1$ such that
$Q_{r;i_1\cdots i_\ell}$ is anti-symmetric with respect to the permutation
of $i_1,\cdots,i_\ell$.
We introduce the lexicographical order from the left on the set
\bea
&&
J_{2n,\ell}=
\{\, (i_1,\cdots,i_\ell)\,\vert\, 
0\leq i_1<\cdots< i_\ell\leq 2n-1\,\}.
\non
\ena
Let us define the $\binom{2n}{\ell}$ by $\binom{2n}{\ell}$ matrix 
$X^{(2n,\ell)}$
by
\bea
&&
X^{(2n,\ell)}=(X^{(2n,\ell)}_{r;i_1\cdots i_\ell}),
\quad
1\leq r\leq \binom{2n}{\ell},
\,
(i_1,\cdots,i_\ell)\in J_{2n,\ell},
\label{matx}
\\
&&
X^{(2n,\ell)}_{r;i_1\cdots i_\ell}=
Q_{r;i_1\cdots i_\ell}.
\non
\ena

\begin{prop}\label{degX}
As an element of $R_{2n}$, $\det X^{(2n,\ell)}$ is homogeneous of degree 
$\binom{2n}{2}\big(\binom{2n-1}{\ell-1}+\binom{2n-2}{\ell-1}\big)$.
\end{prop}
\vskip2mm

\noindent
\pf By the definition of $\deg_1$, if $Q_r=v_I\wedge w_J\wedge \xi_K$,
\bea
&&
\deg_1\, Q_{r;s_1\cdots s_\ell}=d_{I,J,K}+s_1+\cdots+s_\ell.
\non
\ena
Thus, as in the proof of Lemma \ref{lem2}, $\det X^{(2n,\ell)}$ is
homogeneous of degree
\bea
&&
d^{(2n,\ell)}:=\sum_{I,J,K}d_{I,J,K}+
\sum_{0\leq s_1<\cdots<s_\ell\leq 2n-1}(s_1+\cdots+s_\ell).
\non
\ena
The proposition follows from the following lemma.

\begin{lemma}\label{lem-d2nl}
\bea
&&
d^{(2n,\ell)}=\binom{2n}{2}
\Biggl(\binom{2n-1}{\ell-1}+\binom{2n-2}{\ell-1}\Biggr).
\label{d2nl}
\ena
\end{lemma}
\vskip2mm
\noindent
\pf Let us set
\bea
&&
a_{n,\ell}:=\qbc{2n}{\ell}.
\non
\ena
By Proposition \ref{q-tetra-id}, (\ref{ch2nl-1}), (\ref{ch2nl-2})
\bea
&&
a_{n,\ell}=\ch_{2n,\ell}.
\non
\ena
It follows that
\bea
&&
d^{(2n,\ell)}=
\frac{d a_{n,\ell}}{d q}|_{q=1}
+
\sum_{0\leq s_1<\cdots<s_\ell\leq 2n-1}(s_1+\cdots+s_\ell).
\label{d2nl-1}
\ena
We denote the first term of the right hand side of (\ref{d2nl-1}) by
$A_{n,\ell}$ and the second term by $b_{n,\ell}$.
The sum $b_{n,\ell}$ is expressed in terms of $A_{n,\ell}$ as follows.
Notice that
\bea
&&
\sum_{0\leq s_1<\cdots<s_\ell\leq 2n-1}q^{s_1+\cdots+s_\ell}
=q^{\frac{1}{2}\ell(\ell-1)}a_{n,\ell}.
\non
\ena
Then
\bea
&&
b_{n,\ell}=\frac{d}{dq}(q^{\frac{1}{2}\ell(\ell-1)}a_{n,\ell})|_{q=1}
=\binom{\ell}{2}\binom{2n}{\ell}+A_{n,\ell},
\non
\ena
and
\bea
&&
d^{(2n,\ell)}=2A_{n,\ell}+\binom{\ell}{2}\binom{2n}{\ell}.
\non
\ena
Let us define $c_{n,\ell}$ by
\bea
&&
c_{n,\ell}=
\binom{2n}{2}\Biggl(\binom{2n-1}{\ell-1}+\binom{2n-2}{\ell-1}\Biggr)
-
\binom{\ell}{2}\binom{2n}{\ell},
\label{def-c}
\ena
Then the lemma is equivalent to 
\bea
&&
c_{n,\ell}=2A_{n,\ell}.
\label{substep}
\ena
Let us prove (\ref{substep}) by the induction on $n$.
For $n=1$, (\ref{substep}) is easily verified.
Let us prove that $c_{n,\ell}$ and $2A_{n,\ell}$ satisfy
the same recursion relation.

Notice that $a_{n,\ell}$ here is the same as $a_{n,\ell}$ in the 
proof of Proposition \ref{q-tetra-id}.
By differentiating the recursion relation (\ref{a-rec}) of $a_{n,\ell}$
one obtains the recursion relation for $2A_{n,\ell}$:
\bea
2A_{n,\ell}&=&2A_{n-1,\ell}+2(2A_{n-1,\ell-1})+2A_{n-1,\ell-2}
\non
\\
&&
+2(4n-\ell-2)\binom{2n-2}{\ell-1}+2(2n-\ell)\binom{2n-2}{\ell-2}.
\label{rec-A}
\ena

We prove the same recursion relation for $c_{n,\ell}$.
Set
\bea
&&
c^{(1)}_{n,\ell}=\binom{2n}{2}
\Biggl(
\binom{2n-1}{\ell-1}+\binom{2n-2}{\ell-1}
\Biggr),
\quad
c^{(2)}_{n,\ell}=\binom{\ell}{2}\binom{2n}{\ell}.
\non
\ena
Using 
\bea
&&
\binom{n}{m}=\binom{n-2}{m}+2\binom{n-2}{m-1}+\binom{n-2}{m-2},
\non
\ena
in each binomial coefficient in $c^{(1)}_{n,\ell}$ and $c^{(2)}_{n,\ell}$, 
we get
\bea
&&
c^{(1)}_{n,\ell}=
c^{(1)}_{n-1,\ell}+2c^{(1)}_{n-1,\ell-1}+c^{(1)}_{n-1,\ell-2}
+
(4n-3)\Biggl(\binom{2n-1}{\ell-1}+\binom{2n-2}{\ell-2}\Biggr),
\non
\\
&&
c^{(2)}_{n,\ell}=
c^{(2)}_{n-1,\ell}+2c^{(2)}_{n-1,\ell-1}+c^{(2)}_{n-1,\ell-2}
+2(\ell-1)\binom{2n-2}{\ell-1}+(2\ell-3)\binom{2n-2}{\ell-2}.
\non
\ena
Taking the difference of these equations and using the induction 
hypothesis we get the desired recursion relation.
Q.E.D.

\begin{prop}\label{divisible}
The determinant $\det X^{(2n,\ell)}$ is divisible by
$(\Delta_{2n}^{+})^{\binom{2n-1}{\ell-1}+\binom{2n-2}{\ell-1}}$.
\end{prop}
\vskip2mm

\noindent
\pf In $\det X^{(2n,\ell)}$, for each $1\leq j_2<\cdots<j_\ell\leq 2n-1$,
we add the $(2j,j_2,\cdots,j_\ell)$-th
column times $x_{2n-1}^{-2j}$ or $-x_{2n-1}^{-2j}$ to 
the $(0,j_2,\cdots,j_\ell)$-th 
column successively so that, in the resulting determinant, 
the $(r;0,j_2,\cdots,j_\ell)$-th component becomes
\bea
&&
Q'_{r;0,j_2\cdots j_\ell}=
\sum_{j=0}^{n-1} Q_{r;2j,j_2\cdots j_\ell}x_{2n-1}^{-2j}.
\quad
\non
\ena
Further, for each $2\leq j_2'<\cdots<j_\ell'\leq 2n-1$, we add  
the $(2j-1,j_2',\cdots,j_\ell')$-th column times 
$x_{2n-1}^{-(2j-1)}$ or $-x_{2n-1}^{-(2j-1)}$ to 
the $(1,j_2',\cdots,j_\ell')$-th 
column successively so that the $(r;1,j_2',\cdots,j_\ell')$-th 
component becomes
\bea
&&
Q'_{r;1,j_2'\cdots j_\ell'}=
\sum_{j=1}^{n} Q_{r;2j-1,j_2'\cdots j_\ell'}x_{2n-1}^{-(2j-1)}.
\quad
\non
\ena
Since $\{Q_{r;i_1\cdots i_\ell}\}$ satisfy (\ref{eq-P-comp}), we have
\bea
&&
\overline{Q'_{r;0,j_2\cdots j_\ell}}=
\overline{Q'_{r;1,j_2'\cdots j_\ell'}}=0.
\non
\ena
It follows that $Q'_{r;0,j_2\cdots j_\ell}$ and 
$Q'_{r;1,j_2'\cdots j_\ell'}$ are divisible
by $x_{2n-1}+x_{2n}$. 
Thus $\det X^{(2n,\ell)}$ is divisible by
$$
(x_{2n-1}+x_{2n})^{\binom{2n-1}{\ell-1}+\binom{2n-2}{\ell-1}}.
$$
Since $\det X^{(2n,\ell)}$ is a symmetric polynomial of $x_1,...,x_{2n}$,
it is divisible by
$$
(\Delta_{2n}^{+})^{\binom{2n-1}{\ell-1}+\binom{2n-2}{\ell-1}}.
$$
Q.E.D.

\begin{cor}\label{X=D}
There is a constant $c$ such that
\bea
&&
\det X^{(2n,\ell)}=
c\,(\Delta_{2n}^{+})^{\binom{2n-1}{\ell-1}+\binom{2n-2}{\ell-1}}.
\non
\ena
\end{cor}
\vskip2mm

\noindent
Proof of Theorem \ref{conj-basis}.
\par
\noindent
It is sufficient to prove that any $P\in U_{2n,\ell}$ can be written as
a linear combination of $Q_r$'s with the coefficients in $R_{2n}$.
Let us write
\bea
&&
P=
\sum_{0\leq j_1<\cdots<j_\ell\leq 2n-1}
P_{j_1\cdots j_\ell}
X^{j_1}\wedge\cdots\wedge X^{j_\ell}.
\non
\ena
We solve the following linear equations for 
$S_r$, $r=1,\cdots,\binom{2n}{\ell}$:
\bea
&&
P=\sum_{r}S_rQ_r,
\non
\ena
which, in the matrix form,
\bea
&&
(P_{j_1,\cdots,j_\ell})=(S_r)X^{(2n,\ell)},
\non
\ena
where $(P_{j_1,\cdots,j_\ell})$ and $(S_r)$ are both row vectors.
Since $\{Q_r\}$ are linearly independent over $R_{2n}$ by 
Theorem \ref{conj-indep}, $c\neq 0$ in Corollary \ref{X=D}. 
Then by Cramer's formula
\bea
&&
S_r=c^{-1}
(\Delta_{2n}^{+})^{-\binom{2n-1}{\ell-1}-\binom{2n-2}{\ell-1}} 
\det X^{(2n,\ell)}_r,
\non
\ena
where the matrix $\det X^{(2n,\ell)}_r$ is defined by replacing the $r$-th row
in $X^{(2n,\ell)}$ by $(P_{j_1,\cdots,j_\ell})$.
In a similar way to the proof of Proposition \ref{divisible},
$\det X^{(2n,\ell)}_r$ is proved to be divisible by
\bea
&&
(\Delta_{2n}^{+})^{\binom{2n-1}{\ell-1}+\binom{2n-2}{\ell-1}}.
\non
\ena
Thus $S_r\in R_{2n}$ as desired.
Q.E.D.

\section{Equivalent Conditions for linear independence}
Let $\Gamma^{(2n)}$ be the $2n$-dimensional vector space with the basis
$\alpha_i,\beta_j$, $i,j=1,\cdots,n$:
\bea
&&
\Gamma^{(2n)}=\oplus_{i=1}^n\mathbb{C}\alpha_i\oplus_{j=1}^n\mathbb{C}\beta_j.
\non
\ena
Define $\omega_k\in \wedge^2 \Gamma^{(2n)}$ by
\bea
&&
\omega_k=\sum_{r=1}^n \alpha_{r-k+1}\wedge \beta_r.
\quad
1\leq r\leq n,
\label{omega}
\ena
Here and in the following the indices of $\alpha_i$, $\beta_j$, $\omega_k$
are extended to integers and they are read by modulo $n$ in the representative
$\{1,\cdots,n\}$.

In this section we set 
\bea
&&
e_{k}=e_{k}^{(2n)},
\quad
P_{r,s}=P_{r,s}^{(2n)},
\quad
v_i=v_i^{(2n)},
\quad
w_j=w_j^{(2n)},
\quad
\xi_k=\xi_k^{(2n)},
\non
\ena
since $n$ always remains same.

\begin{theorem}\label{equiv-cond}
The following conditions are equivalent:
\begin{itemize}
\item[(i)]
$\{\alpha_I\wedge \beta_J\wedge \omega_K\}$ is linearly independent 
over $\mathbb{C}$,
\item[(ii)]
$\{\alpha_I\wedge \beta_J\wedge \omega_K\}$ linearly generate 
$\wedge^\ell \Gamma^{(2n)}$,
\item[(iii)]
$\{\alpha_I\wedge \beta_J\wedge \omega_K\}$ is a basis of
$\wedge^\ell \Gamma^{(2n)}$,
\item[(iv)]
$\{v_I\wedge w_J\wedge \xi_K\}$ is linearly independent over $R_{2n}$,
\end{itemize}
where $I,J,K$ satisfy (\ref{indexI})-(\ref{sum=l}).
\end{theorem}
\vskip2mm

\noindent
\pf By (\ref{dim-correct}), it is obvious that
(i), (ii) and (iii) are equivalent.
Let us prove the equivalence of (i) and (iv).
To this end we have to prepare several lemmas and propositions.

\begin{prop}\label{prop71}
The following expressions for $\xi_k$ hold:
\bea
&&
\xi_k=
\Big(
-\sum_{i=1}^{k-1}\sum_{r=n-k+1}^{n-k+i}
+\sum_{i=k}^{n}\sum_{r=i-k}^{n-k}
\Big)
e_{2(k-i+r)}X^{2r+1}\wedge v_i,
\label{xiexp-1}
\\
&&
\xi_k=\Big(
-\sum_{i=1}^{k-1}\sum_{r=n-k+1}^{n-k+i}
+\sum_{i=k}^{n}\sum_{r=i-k}^{n-k}
\Big)
e_{2(k-i+r)}w_i\wedge X^{2r}.
\label{xiexp-2}
\ena
\end{prop}

The proof of this proposition is given in the end of this section.

We consider the specialization
\bea
&&
e_1=-e_{2n}=1,
\quad
e_j=0,
\quad
j\neq1,2n.
\label{special}
\ena

\begin{lemma}\label{lem71}
At (\ref{special}) we have
$$
P_{r,s}=
\left\{
\begin{array}{rl}
1,&\quad \text{$1\leq s\leq n$ and $s=n-r+2$ mod. $n$}\\
0,&\quad \text{otherwise}.
\end{array}
\right.
$$\end{lemma}
\vskip2mm

\noindent
\pf Let us prove the statement by the induction on $r$.
For $r=1$, 
\bea
&&
P_{1,s}=e_{2s-1},
\non
\ena
and the statement is obvious.
For $r\geq 2$, Substituting (\ref{special}) into (\ref{def-P}) and
using the hypothesis of the induction we have
$$
P_{r,n-r+2}=P_{r-1,n-(r-1)+2}-e_{2(n-r+2)}P_{r-1,1}
=\left\{
\begin{array}{rl}
P_{1,n+1}-e_{2n}P_{1,1}=1,&\quad r=2\\
1,&\quad r\geq 3.
\end{array}
\right.
$$
For $s\neq n-r+2$,
\bea
&&
P_{r,s}=P_{r-1,s+1}-e_{2s}P_{r-1,1}=0,
\non
\ena
by the induction hypothesis.
Q.E.D.

For any integer $i$ we set
\bea
&&
<i>=X^j,
\quad
0\leq j\leq 2n-1,
\quad
i\equiv j \text{ mod. $2n$.}
\non
\ena

\begin{cor}\label{cor71}
At (\ref{special}) we have
\bea
&&
v_i=<2(n-i+1)>,
\quad
w_i=<2(n-i+1)+1>,
\quad
1\leq i\leq n.
\non
\ena
\end{cor}

The following lemma follows from Proposition \ref{prop71}
and Corollary \ref{cor71}.

\begin{lemma}\label{lemma72}
At (\ref{special}),
\bea
&&
\xi_k=\sum_{r=1}^n<2(n+r-k)+1>\wedge v_r=\sum_{r=1}^n w_{k-r+1}\wedge v_r.
\non
\ena
\end{lemma}

At (\ref{special}) we set
\bea
&&
\alpha_i=v_{2-i},
\quad
\beta_j=w_j.
\non
\ena
Then we have
\bea
&&
\omega_k=\sum_{r=1}^n \alpha_{r-k+1}\wedge \beta_r=-\xi_k.
\non
\ena
Thus at (\ref{special}) 
\bea
&&
\{v_I\wedge w_J \wedge \xi_K\}=\{\pm \alpha_I\wedge \beta_J \wedge \omega_K\},
\non
\ena
where $\pm$ means that $+$ or $-$ is determined to each $(I,J,K)$.
Consequently, if $\{\alpha_I\wedge \beta_J \wedge \omega_K\}$ are linearly
independent, for the matrix $X^{(2n,\ell)}$ defined in (\ref{matx}),
\bea
&&
\det X^{(2n,\ell)}\neq 0,
\non
\ena
at (\ref{special}) and $\{v_I\wedge w_J \wedge \xi_K\}$ are linearly
independent over $R_{2n}$. Therefore (iv)$\Rightarrow$(i) in 
Theorem \ref{equiv-cond} is proved. 

Let us prove the converse. We have to use the following lemma
which can be easily verified.

\begin{lemma}\label{lem73}
At (\ref{special}), we have
\bea
&&
\Delta^{+}_{2n}=(-1)^{\frac{1}{2}n(n+1)
(\binom{2n-1}{\ell-1}+\binom{2n-2}{\ell-1})}.
\non
\ena
\end{lemma}
\vskip2mm

Now, if $\{v_I\wedge w_J \wedge \xi_K\}$ are linearly independent over 
$R_{2n}$, $c\neq 0$ in Corollary \ref{X=D}.
Then, at (\ref{special}), we get
\bea
&&
\det X^{(2n,\ell)}=c(-1)^{\frac{1}{2}n(n+1)
(\binom{2n-1}{\ell-1}+\binom{2n-2}{\ell-1})}\neq 0.
\non
\ena
Thus $\{\alpha_I\wedge \beta_J \wedge \omega_K\}$ are linearly
independent.
Q.E.D.
\vskip5mm

\noindent
{\it Proof of Proposition \ref{prop71}}.
\par
\noindent
Let us first prove (\ref{xiexp-1}).
In this proof we extend the index $r$ of $P_{r,s}$ to $r\in \mathbb{Z}$
by defining
\bea
&&
P_{r,s}=e_{2(s+r)-3},
\quad
r\leq 0.
\non
\ena
For $r\geq n+1$, $P_{r,s}$ is defined by the recursion relation (\ref{def-P}).
Then the recursion relation (\ref{def-P}) holds for all $r,s\in \mathbb{Z}$.

The following lemma is obtained by a simple calculation 
using the definitions (\ref{def-vw}) and (\ref{def-xi}).

\begin{lemma}\label{lem74}
\bea
&&
(X_1+X_2)\xi_k(X_1,X_2)=
\sum_{i=0}^n\sum_{j=0}^n(e_{2j}P_{k,i}-e_{2i}P_{k,j})X_1^{2i}X_2^{2j}.
\non
\ena
\end{lemma}

Let us denote the right hand side of (\ref{xiexp-1}) by $\eta_k$
and expand it as
\bea
&&
\eta_k=\sum_{i,j=0}^{n-1} a^{(k)}_{ij}(X_1^{2i+1}X_2^{2j}-X_2^{2i+1}X_1^{2j}).
\non
\ena
We set $a^{(k)}_{ij}=0$ unless $0\leq i,j\leq n-1$.
Then
\bea
&&
(X_1+X_2)\eta_k=
\sum_{i,j}(a^{(k)}_{i-1,j}-a^{(k)}_{j-1,i})X_1^{2i}X_2^{2j}
+\sum_{i,j}(a^{(k)}_{ij}-a^{(k)}_{ji})X_1^{2i+1}X_2^{2j+1}.
\non
\ena
Comparing this with Lemma \ref{lem74} we see that (\ref{xiexp-1}) is equivalent
to the following equations:
\bea
&&
a^{(k)}_{ij}=a^{(k)}_{ji},
\label{eq-72}
\\
&&
a^{(k)}_{i-1,j}-a^{(k)}_{j-1,i}=e_{2j}P_{k,i}-e_{2i}P_{k,j}.
\label{eq-73}
\ena
In order to prove these equations we calculate $a^{(k)}_{ij}$
more explicitly.

\begin{lemma}\label{lem75}
We have
$$
a^{(k)}_{ij}=
\left\{
\begin{array}{rl}
\sum_{r=0}^i e_{2(i-r)}P_{k+r,j+1},&
\quad
0\leq i\leq n-k
\\
-\sum_{r=1}^{n-i} e_{2(i+r)}P_{k-r,j+1},&
\quad
n-k+1\leq i\leq n-1.
\end{array}
\right.
$$
\end{lemma}
\vskip2mm

\noindent
\pf If $0\leq i\leq n-k$, $X^{2i+1}\wedge X^{2j}$ appears only from the
second sum of the right hand side of (\ref{xiexp-1}). Its coefficient is
calculated as
\bea
&&
a^{(k)}_{ij}=
\sum_{r=k}^{k+i}e_{2(k-r+i)}P_{r,j+1}
=
\sum_{r=0}^{i}e_{2(i-r)}P_{k+r,j+1}.
\non
\ena
Similarly, if $n-k+1\leq i\leq n-1$, $X^{2i+1}\wedge X^{2j}$ 
appears only from the
first sum term and its coefficient is given by
\bea
&&
a^{(k)}_{ij}=
-\sum_{r=k-n+i}^{k-1}e_{2(k-r+i)}P_{r,j+1}
=
-\sum_{r=1}^{n-i}e_{2(i+r)}P_{k-r,j+1}.
\non
\ena
Q.E.D.

Let us rewrite $a^{(k)}_{ij}$ more.

\begin{lemma}\label{lem76}
$$
a^{(k)}_{ij}=
\left\{
\begin{array}{ll}
\big(
\sum_{r+s=i+j+2,\,r\leq i}-\sum_{r+s=i+j+2,\,r\geq j+1}
\big)
e_{2r}P_{k-1,s},
&
\quad 
0\leq i\leq n-k
\\
\big(
\sum_{r+s=i+j+1,\,r\leq j}-\sum_{r+s=i+j+1,\,r\geq i+1}
\big)
e_{2r}P_{k,s},
&
\quad 
n-k+1\leq i\leq n-1.
\end{array}
\right.
$$
\end{lemma}
\vskip2mm

\noindent
\pf
First we suppose that $0\leq i\leq n-k$.
Using Lemma \ref{lem75} we have
\bea
a^{(k)}_{ij}&=&\sum_{r=0}^ie_{2(i-r)}P_{k+r,j+1}
\non
\\
&=&
\sum_{r=0}^ie_{2(i-r)}(P_{k+r-1,j+2}-e_{2(j+1)}P_{k+r-1,1})
\non
\\
&=&
\sum_{r=1}^ie_{2(i-r)}P_{k+r-1,j+2}+e_{2i}P_{k-1,j+2}
-e_{2(j+1)}\sum_{r=0}^ie_{2(i-r)}P_{k+r-1,1}.
\non
\ena
By Lemma \ref{lem75} we have, for the first term,
\bea
&&
\sum_{r=1}^ie_{2(i-r)}P_{k+r-1,j+2}=a^{(k)}_{i-1,j+1}.
\non
\ena
As to the third term, we see
\bea
\sum_{r=0}^ie_{2(i-r)}P_{k+r-1,1}
&=&
P_{k+i-1,1}+e_2P_{k+i-2,1}+\cdots+e_{2i}P_{k-1,1}
\non
\\
&=&
(P_{k+i-2,2}-e_2P_{k+i-2,1})+e_2P_{k+i-2,1}+\cdots+e_{2i}P_{k-1,1}
\non
\\
&=&
P_{k+i-2,2}+e_4P_{k+i-3,1}+\cdots+e_{2i}P_{k-1,1}
\non
\\
&\vdots&
\non
\\
&=&
P_{k-1,i+1}.
\non
\ena
Thus we have
\bea
a^{(k)}_{ij}
&=&a^{(k)}_{i-1,j+1}+e_{2i}P_{k-1,j+2}-e_{2(j+1)}P_{k-1,i+1}
\non
\\
&=&
a^{(k)}_{0,i+j}+
\sum_{r=0}^{i-1}(e_{2(i-r)}P_{k-1,j+2+r}
-e_{2(j+1+r)}P_{k-1,i+1-r}).
\label{eq-74}
\ena
Notice that,
\bea
&&
a^{(k)}_{0,i+j}=P_{k,i+j+1}=P_{k-1,i+j+2}-e_{2(i+j+1)}P_{k-1,1},
\non
\ena
which is precisely the $r=i$ term in the second term of (\ref{eq-74}).
Therefore
\bea
&&
a^{(k)}_{ij}=
\sum_{r=0}^{i}
(e_{2(i-r)}P_{k-1,j+2+r}-e_{2(j+1+r)}P_{k-1,i+1-r}),
\non
\ena
which proves the first equation of the lemma.

Next we consider the case $n-k+1\leq i\leq n-1$.
Here we use the recursion relation of $P_{r,s}$ in the form
\bea
&&
P_{r,s}=P_{r+1,s-1}+e_{2(s-1)}P_{r,1}.
\label{eq-75}
\ena
Using (\ref{eq-75}) we have
\bea
a^{(k)}_{ij}
&=&
-\sum_{r=1}^{n-i}e_{2(i+r)}P_{k-r,j+1}
\non
\\
&=&
-\sum_{r=1}^{n-i}e_{2(i+r)}(P_{k-r+1,j}+e_{2j}P_{k-r,1})
\non
\\
&=&
-\sum_{r=2}^{n-i}e_{2(i+r)}P_{k-r+1,j}
-e_{2(i+1)}P_{k,j}
-e_{2j}\sum_{r=1}^{n-i}e_{2(i+r)}P_{k-r,1},
\non
\ena
where we have
\bea
&&
-\sum_{r=2}^{n-i}e_{2(i+r)}P_{k-r+1,j}=a^{(k)}_{i+1,j-1},
\non
\ena
and
\bea
\sum_{r=1}^{n-i}e_{2(i+r)}P_{k-r,1}
&=&
P_{k,i+1}+\sum_{r=1}^{n-i}e_{2(i+r)}P_{k-r,1}-P_{k,i+1}
\non
\\
&=&
P_{k-n+i,n+1}-P_{k,i+1}
\non
\\
&=&
-P_{k,i+1},
\non
\ena
where we use Lemma \ref{lem0}.
Therefore
\bea
a^{(k)}_{ij}
&=&
a^{(k)}_{i+1,j-1}+e_{2j}P_{k,i+1}-e_{2(i+1)}P_{k,j}
\non
\\
&=&
a^{(k)}_{n-1,i+j+1-n}
+\sum_{r=0}^{n-i-2}(e_{2(j-r)}P_{k,i+1+r}-e_{2(i+1+r)}P_{k,j-r}).
\label{eq-76}
\ena
Here
\bea
&&
a^{(k)}_{n-1,i+j+1-n}=-e_{2n}P_{k-1,i+j+2-n},
\non
\ena
which is the $r=n-i-1$ term in the second term of (\ref{eq-76}).
Thus we have 
\bea
&&
a^{(k)}_{ij}=
\sum_{r=0}^{n-i-1}(e_{2(j-r)}P_{k,i+1+r}-e_{2(i+1+r)}P_{k,j-r}),
\non
\ena
which completes the proof of Lemma \ref{lem76}.
Q.E.D.
\vskip5mm

\noindent
{\it Proof of (\ref{eq-72})}.
\par
\noindent
We can assume $i<j$. The proof is divided into three cases.
\vskip2mm

\noindent
(i). the case $0\leq i<j\leq n-k$:
We use the first equation in Lemma \ref{lem76} and notations like
\bea
&&
(r\leq i)_1=\sum_{r+s=i+j+2,\,r\leq i}e_{2r}P_{k-1,s}.
\label{eq-77}
\ena
Then we have
\bea
&&
a^{(k)}_{ij}=(r\leq i)_1-(r\geq j+1)_1,
\non
\ena
and
\bea
&&
a^{(k)}_{ji}=(r\leq j)_1-(r\geq i+1)_1=(r\leq i)_1-(r\geq j+1)_1=a^{(k)}_{ij}.
\non
\ena
\vskip2mm

\noindent
(ii). the case $n-k+1\leq i<j\leq n-1$:
In this case we use notations like
\bea
(r\leq j)_2=\sum_{r+s=i+j+1,\,r\leq j}e_{2r}P_{k,s}.
\label{eq-78}
\ena
Then
\bea
&&
a^{(k)}_{ij}=(r\leq j)_2-(r\geq i+1)_2=(r\leq i)_2-(r\geq j+1)_2=a^{(k)}_{ji}.
\non
\ena
\vskip2mm

\noindent
(iii). the case $0\leq i\leq n-k<j\leq n-1$:
In this case we have
\bea
a^{(k)}_{ji}
&=&
(\sum_{r+s=i+j+1,\,r\leq i}-\sum_{r+s=i+j+1,\,r\geq j+1})e_{2r}P_{k,s}
\non
\\
&=&
(\sum_{r+s=i+j+1,\,r\leq i}-\sum_{r+s=i+j+1,\,r\geq j+1})e_{2r}
(P_{k-1,s+1}-e_{2s}P_{k-1,1})
\non
\\
&=&
(\sum_{r+s=i+j+2,\,r\leq i,\,s\geq 1}-\sum_{r+s=i+j+2,\,r\geq j+1,\,s\geq 1})
e_{2r}(P_{k-1,s}-e_{2(s-1)}P_{k-1,1}).
\non
\ena
Notice that, if $s=0$,
\bea
&&
P_{k-1,s}-e_{2(s-1)}P_{k-1,1}=0.
\non
\ena
Then
\bea
&&
a^{(k)}_{ji}=
a^{(k)}_{ij}-
(\sum_{r+s=i+j+1,\,r\leq i}-\sum_{r+s=i+j+1,\,r\geq j+1})e_{2r}e_{2s}P_{k-1,1}.
\non
\ena
The second term in the right hand side vanishes because $r\leq i$ is equivalent to $s\geq j+1$ under the condition $r+s=i+j+1$. 
Thus (\ref{eq-72}) is proved.
\vskip3mm

\noindent
{\it Proof of (\ref{eq-73})}.
\par
\noindent
Since both hand sides of (\ref{eq-73}) are anti-symmetric with respect to
$i$ and $j$, we can assume $i<j$.

\noindent
(i). the case $0\leq i\leq n-k$:
We use the notation (\ref{eq-77}) in which $r+s=i+j+2$ 
is changed to $r+s=i+j+1$.
Then we have
\bea
&&
a^{(k)}_{i-1,j}=(r\leq i-1)_1-(r\geq j+1)_1,
\non
\\
&&
a^{(k)}_{j-1,i}=a^{(k)}_{i,j-1}=(r\leq i)_1-(r\geq j)_1.
\non
\ena
Thus
\bea
a^{(k)}_{i-1,j}-a^{(k)}_{j-1,i}
&=&
-(r=i)_1+(r=j)_1
\non
\\
&=&
-e_{2i}P_{k-1,j+1}+e_{2j}P_{k-1,i+1}
\non
\\
&=&
-e_{2i}P_{k,j}+e_{2j}P_{k,i},
\non
\ena
where to derive the last equation we have used the recursion relation of 
$P_{r,s}$.

\noindent
(ii). the case $n-k+1\leq i<j\leq n-1$:
We use the notation (\ref{eq-78}) in which $r+s=i+j+1$ is changed to $r+s=i+j$.
Then
\bea
&&
a^{(k)}_{i-1,j}=a^{(k)}_{j,i-1}=(r\leq i-1)_2-(r\geq j+1)_2,
\non
\\
&&
a^{(k)}_{j-1,i}=(r\leq i)_2-(r\geq j)_2.
\non
\ena
Therefore
\bea
a^{(k)}_{i-1,j}-a^{(k)}_{j-1,i}
&=&
-(r=i)_2+(r=j)_2
\non
\\
&=&
-e_{2i}P_{k,j}+e_{2j}P_{k,i},
\non
\ena
which completes the proof of (\ref{eq-73}).
\vskip2mm

Thus (\ref{xiexp-1}) is proved.

The expression (\ref{xiexp-2}) follows from (\ref{xiexp-1}) in the following
manner.

Let $\eta_k'$ denote the right hand side of (\ref{xiexp-2}).
Notice that $w_i(X)$ is obtained from $v_i(X)$ by the replacement
\bea
&&
X^{2r}\ra X^{2r+1},
\quad
0\leq r\leq n-1.
\non
\ena
Then it is obvious that $\eta_k'$ is obtained from $-\eta_k$
by the replacement
\bea
&&
X^{2r}\ra X^{2r+1},
\quad
X^{2r+1}\ra X^{2r},
\quad
0\leq r\leq n-1.
\label{eq-79}
\ena
Since
\bea
&&
\eta_k=\sum a^{(k)}_{ij}X^{2i+1} \wedge X^{2j},
\non
\ena
we have
\bea
\eta_k'&=& -\sum a^{(k)}_{ij}X^{2i} \wedge X^{2j+1}
\non
\\
&=&
\sum a^{(k)}_{ji}X^{2i+1} \wedge X^{2j}
\non
\\
&=&
\sum a^{(k)}_{ij}X^{2i+1} \wedge X^{2j}=\eta_k,
\non
\ena
which proves (\ref{xiexp-2}). Thus the proof of Proposition \ref{prop71}
is completed.
Q.E.D.

\section{Linear independence}\label{pf-indep}
In this section we prove 

\begin{theorem}\label{independence}
The elements $\{\alpha_I\wedge \beta_J\wedge \omega_K\}$ linearly generate
$\wedge^\ell \Gamma^{(2n)}$, where $I,J,K$ satisfy 
(\ref{indexI})-(\ref{sum=l}).
\end{theorem}
\vskip1.5mm

As a consequence of Theorem \ref{independence} and \ref{equiv-cond}, 
Theorem \ref{conj-indep} is proved.
\vskip2mm

Let $B_\ell$ be the linear span of $\{\alpha_I\wedge \beta_J\wedge \omega_K\}$
in $\wedge^\ell \Gamma^{(2n)}$. We use the symbol $\equiv$ for the equality
modulo $B_\ell$. For the sake of simplicity, in the following,
we omit the symbol $\wedge$ of the exterior product and simply write 
the product for the exterior product.

The conditions (\ref{indexI})-(\ref{sum=l}) for the product $\alpha_I\beta_J$,
$I=(i_1,\cdots,i_{\ell_1})$, $J=(j_1,\cdots,j_{\ell_2})$ are
\bea
&&
1\leq i_1\leq \cdots \leq i_{\ell_1}\leq n,
\quad
1\leq j_1\leq \cdots \leq j_{\ell_2}\leq n-\ell_1,
\quad
\ell_1+\ell_2=\ell.
\non
\ena
The elements $\alpha_I\beta_J$ which does not satisfy 
(\ref{indexI})-(\ref{sum=l}) are classified by the number of $k$
of $\beta_j$ such that $j>n-\ell_1$. We write such $j$ in the form
$j=n-r$ with $0\leq r\leq \ell_1-1$ for later convenience.
With these notations, in order to prove Theorem \ref{independence},
it is sufficient to prove  
\bea
&&
\alpha_{i_{\ell_1-1}}\cdots\alpha_{i_{0}}
\beta_{j_1}\cdots\beta_{j_{\ell_2-k-1}}
\beta_{n-r_k}\cdots\beta_{n-r_0}
\equiv0,
\label{eq9-1}
\\
&&
1\leq i_{\ell_1-1}<\cdots<i_0\leq n,
\quad
1\leq j_1<\cdots<j_{\ell_2-k-1}\leq n-\ell_1,
\non
\\
&&
0\leq r_0<\cdots<r_k\leq \ell_1-1,
\quad
0\leq k\leq \ell_1-1,
\non
\ena
here we changed the way of numbering the indices of $\alpha$ for the sake
of later convenience.
We set
\bea
&&
J'=(j_1,\cdots,j_{\ell_2-k-1}).
\non
\ena
Let us consider the element
\bea
&&
\alpha_{i_{\ell_1-1}}\cdots
\widehat{\alpha_{i_{r_k}}}\cdots
\widehat{\alpha_{i_{r_0}}}\cdots
\alpha_{i_{0}}
\beta_{J'}
\omega_{n-r_k-i_{r_k}+1}
\cdots
\omega_{n-r_0-i_{r_0}+1},
\label{eq9-2}
\ena
where $\widehat{\alpha_{i_{r_k}}}$ means that $\alpha_{i_{r_k}}$ is 
removed.
Since
\bea
&&
i_s\geq \ell_1-s,
\quad
0\leq s \leq \ell_1-1,
\non
\ena
we have
\bea
&&
n+1-r_s-i_{r_s}\leq n+1-\ell_1,
\quad
0\leq s \leq k.
\non
\ena
Notice that, in (\ref{eq9-2}), the number of $\alpha$ is $\ell_1-k-1$ and
the number of $\omega$ is $k+1$. Then the indices of (\ref{eq9-2}) 
satisfy (\ref{indexI})-(\ref{sum=l})
if and only if the indices of $\beta_{j'}$, $\omega_{k'}$ appearing in 
(\ref{eq9-2}) satisfy $j'\leq n-\ell_1$ and $k'\leq n+1-\ell_1$.
Thus 
\bea
&&
\text{(\ref{eq9-2})}\equiv 0.
\label{eq9-3}
\ena
Let us prove (\ref{eq9-1}) by the induction on $k$.
The case $k=0$ is proved in Proposition \ref{prop9-1} below.

Suppose that $k\geq 1$ and (\ref{eq9-1}) is valid until $k-1$.
Then it is possible to expand (\ref{eq9-2}) using (\ref{omega}) as
\bea
\text{(\ref{eq9-2})}
&\equiv&
\alpha_{i_{\ell_1-1}}\cdots
\widehat{\alpha_{i_{r_k}}}\cdots
\widehat{\alpha_{i_{r_0}}}\cdots
\alpha_{i_{0}}
\beta_{J'}
\sum_{p_0,...,p_k=n-\ell_1+1}^n
\alpha_{p_k+r_k+i_{r_k}}\beta_{p_k}
\cdots
\alpha_{p_0+r_0+i_{r_0}}\beta_{p_0}
\non
\\
&\equiv&
\pm
\sum_{p_0,...,p_k=0}^{\ell_1-1}
\alpha_{i_{\ell_1-1}}\cdots
\alpha_{i_{r_k}+r_k-p_k}\cdots
\alpha_{i_{r_0}+r_0-p_0}\cdots
\alpha_{i_{0}}
\beta_{J'}
\beta_{n-p_k}\cdots\beta_{n-p_0}
\non
\\
&\equiv&
\pm
\sum_{0\leq p_0<\cdots<p_k\leq \ell_1-1}
\sum_{\sigma\in S_{k+1}}\sgn\, \sigma \cdot
\alpha_{i_{\ell_1-1}}\cdots
\alpha_{i_{r_k}+r_k-p_{\sigma_k}}\cdots
\alpha_{i_{r_0}+r_0-p_{\sigma_0}}\cdots
\alpha_{i_{0}}
\non
\\
&&
\qquad\quad
\times
\beta_{J'}
\beta_{n-p_k}\cdots\beta_{n-p_0},
\label{eq9-4}
\ena
where the equation is valid for $k=0$ without any assumption and
we use the induction hypothesis for $k\geq 1$.
The sign $\pm$ is determined from  $\ell_2$, $k$, $\{r_s\}$. 
Since it is irrelevant to the following argument, 
we omit its explicit description.

By (\ref{eq9-3}) and (\ref{eq9-4}) we have
\bea
0&\equiv&
\alpha_{i_{\ell_1-1}}\cdots\alpha_{i_{0}}
\beta_{J'}\beta_{n-r_k}\cdots\beta_{n-r_0}
\label{eq9-5'}
\\
&+&
\sum_{\sigma\neq 1}
\sgn\, \sigma \cdot
\alpha_{i_{\ell_1-1}}\cdots
\alpha_{i_{r_k}+r_k-r_{\sigma_k}}\cdots
\alpha_{i_{r_0}+r_0-r_{\sigma_0}}\cdots
\alpha_{i_{0}}
\beta_{J'}
\beta_{n-r_k}\cdots\beta_{n-r_0}
\label{eq9-5}
\\
&+&
\hbox{$\sum^\prime_{0\leq p_0<\cdots<p_k\leq \ell_1-1}
\sum_{\sigma}$}
(\text{the summand of (\ref{eq9-4})}),
\non
\ena
where $\sum^{\prime}$ means that $(p_0,\cdots,p_k)=(r_0,\cdots,r_k)$
is excluded.

\begin{lemma}\label{lem9-1}
If a summand of (\ref{eq9-5}) is proportional to
\bea
&&
\alpha_{i_{\ell_1-1}}\cdots\alpha_{i_{0}}
\beta_{J'}
\beta_{n-r_k}\cdots\beta_{n-r_0},
\label{eq9-6}
\ena
it coincides with (\ref{eq9-6}).
\end{lemma}
\vskip2mm

\noindent
\pf
It is sufficient to prove that, if
\bea
&&
(i_{r_k}+r_k-r_{\sigma_k},\cdots,i_{r_0}+r_0-r_{\sigma_0})
\non
\ena
is a permutation of $(i_{r_k},\cdots,i_{r_0})$, then
\bea
&&
(i_{r_k}+r_k-r_{\sigma_k},\cdots,i_{r_0}+r_0-r_{\sigma_0})
=(i_{r_{\sigma_k}},\cdots,i_{r_{\sigma_0}}).
\label{eq9-7}
\ena
Notice that
\bea
&&
i_{r_k}+R_k \leq \cdots \leq i_{r_1}+R_1\leq i_{r_0},
\quad
R_s=r_s-r_0.
\label{eq9-8}
\ena
We set $j_0=i_{r_0}$ for the sake of simplicity and write 
\bea
&&
(i_{r_k}+R_k,\cdots,i_{r_0}+R_0)=
((j_0-N)^{M_N}\cdots (j_0-1)^{M_1}j_0^{M_0}),
\quad
M_N\neq 0,
\non
\ena
where
\bea
&&
(j_0-s)^{M_s}=j_0-s,\cdots,j_0-s\quad (\text{$M_s$ terms}).
\non
\ena
Then
\bea
&&
(i_{r_k},\cdots,i_{r_0})=(j_0-N-R_k,\cdots,j_0),
\non
\\
&&
(i_{r_k}+r_k-r_{\sigma_k},\cdots,i_{r_0}+r_0-r_{\sigma_0})
=
((j_0-N)^{M_N}\cdots j_0^{M_0})
-
(R_{\sigma_k},\cdots,R_{\sigma_0}).
\label{eq9-9}
\ena
Suppose that (\ref{eq9-9}) is a permutation of $(i_{r_k},\cdots,i_{r_0})$.
Notice the term which coincides with $i_{r_k}$.
If
\bea
&&
j_0-a-R_s=i_{r_k}=j_0-N-R_k,
\non
\ena
for some $0\leq a\leq N$, then obviously we have $a=N$ and $s=k$.
Next consider the term which coincides with $i_{r_{k-1}}$ and so on.
In this way we see that, for each $0\leq s\leq N$, 
\bea
&&
(\sigma_{k-M_N-\cdots-M_s},\cdots,\sigma_{k-M_N-\cdots-M_{s-1}+1})
\non
\ena
is a permutation of 
\bea
&&
(k-M_N-\cdots-M_s,\cdots,k-M_N-\cdots-M_{s-1}+1).
\non
\ena
This means
\bea
&&
((j_0-N)^{M_N}\cdots j_0^{M_0})
-
(R_{\sigma_k},\cdots,R_{\sigma_0})
\non
\\
&&=
(j_0-N-R_{\sigma_k},\cdots,j_0-N-R_{\sigma_{k-M_N+1}},\cdots)
\non
\\
&&=
(i_{r_{\sigma_k}},\cdots,i_{r_{\sigma_0}}).
\non
\ena
Q.E.D.
\vskip2mm

Let $M-1$ be the number of (\ref{eq9-6}) appearing in (\ref{eq9-5}).
Then we can solve the equation (\ref{eq9-5'}) as
\bea
&&
M\alpha_{i_{\ell_1-1}}\cdots\alpha_{i_{0}}
\beta_{J'}\beta_{n-r_k}\cdots\beta_{n-r_0}
\non
\\
&&
\equiv
-
\hbox{$\sum^\prime_{0\leq p_0<\cdots<p_k\leq \ell_1-1}$}
\hbox{$\sum^\prime_{\sigma}$} \sgn\, \sigma \cdot
\alpha_{i_{\ell_1-1}}\cdots
\alpha_{i_{r_k}+r_k-p_{\sigma_k}}\cdots
\alpha_{i_{r_0}+r_0-p_{\sigma_0}}\cdots
\alpha_{i_{0}}
\non
\\
&&
\qquad
\times 
\beta_{J'}
\beta_{n-p_k}\cdots\beta_{n-p_0},
\label{eq9-10}
\ena
where $\sum^\prime$  means that the terms proportional to (\ref{eq9-6})
are excluded from the sum.

In order to simplify the notation we set
\bea
&&
(i_{\ell_1-1},\cdots,i_0|n-r_k,\cdots,n-r_0):=
\alpha_{i_{\ell_1-1}}\cdots\alpha_{i_{0}}
\beta_{J'}
\beta_{n-r_k}\cdots\beta_{n-r_0}.
\label{eq9-11}
\ena
Notice that in the left hand side $J'$ does not appear.
We omit $J'$ because $J'$ remains same in all the arguments below.

Let us consider the map among elements of the form (\ref{eq9-11}) and $0$:
\bea
&&
(i_{\ell_1-1},\cdots,i_0|n-r_k,\cdots,n-r_0)
\non
\\
&&
\mapsto
(i_{\ell_1-1},\cdots, i_{r_k}+r_k-p_{\sigma_k},\cdots,
i_{r_0}+r_0-p_{\sigma_0},\cdots,i_0|
n-p_k,\cdots,n-p_0),
\label{eq9-12}
\ena
where $\{r_s\}$, $\{p_s\}$, $\sigma$, $\{i_s\}$ should satisfy
\bea
&&
0\leq r_0<\cdots< r_k\leq \ell_1-1,
\quad
0\leq p_0<\cdots< p_k\leq \ell_1-1,
\label{eq9-13}
\\
&&
\text{LHS}\neq \pm \text{RHS},
\label{eq9-14}
\\
&&
1\leq i_{\ell_1-1}<\cdots<i_0\leq n.
\non
\ena

For a set of integers $(i_{\ell_1-1},\cdots,i_0)$ we define
\bea
&&
h(i_{\ell_1-1},\cdots,i_0)=
\max\{i_{\ell_1-1},\cdots,i_0\}
-
\min\{i_{\ell_1-1},\cdots,i_0\}.
\non
\ena
For (\ref{eq9-11}) we set
\bea
&&
h(i_{\ell_1-1},\cdots,i_0|n-r_k,\cdots,n-r_0):=
h(i_{\ell_1-1},\cdots,i_0).
\non
\ena

\begin{lemma}\label{lem9-2}
Suppose that the right hand side of (\ref{eq9-12}) is written
as 
\bea
&&
\pm(i'_{\ell_1-1},\cdots,i'_0|n-p_k,\cdots,n-p_0),
\quad
i'_{\ell_1-1}\leq \cdots\leq i_0'.
\non
\ena
Then 
\bea
&&
i_{\ell_1-1}\leq i'_{\ell_1-1},
\quad
i_0'\leq i_0.
\non
\ena
In particular 
\bea
h(i'_{\ell_1-1},\cdots,i'_0|n-p_k,\cdots,n-p_0)
\leq
h(i_{\ell_1-1},\cdots,i_0|n-p_k,\cdots,n-p_0),
\label{eq9-15}
\ena
and the equality holds if and only if
\bea
&&
i_{\ell_1-1}= i'_{\ell_1-1},
\quad
i_0'= i_0.
\non
\ena
\end{lemma}
\vskip2mm

\noindent
\pf
Notice that
\bea
&&
i_s+s\leq i_{s'}+s',
\quad
0\leq s'\leq s\leq \ell_1-1.
\non
\ena
It follows that, for any $0\leq s\leq \ell_1-1$, 
\bea
&&
i_{\ell_1-1}\leq i_{\ell_1-1}+\ell_1-1-p_{\sigma_s}\leq 
i_{r_s}+r_s-p_{\sigma_s}\leq i_0-p_{\sigma_s}\leq i_0,
\non
\ena
which proves the lemma.
Q.E.D.
\vskip2mm

We denote the left hand and the right hand sides of (\ref{eq9-15})
by $h'$ and $h$ respectively.
By Lemma \ref{lem9-2}, the case $h=h'$ occurs only in the following four cases.
\begin{itemize}
\item[(i).] $r_k\neq \ell_1-1$, $r_0\neq 0$.
\item[(ii).] $r_k\neq \ell_1-1$, $r_0=0$ and there exists $k_1$ such that
$i_{r_{k_1}}+r_{k_1}=i_0$.
\item[(iii).] $r_k =\ell_1-1$, $r_0\neq 0$ and there exists $k_1$ such that
$i_{r_{k_1}}+r_{k_1}-(\ell_1-1)=i_{\ell_1-1}$.
\item[(iv).] $r_k =\ell_1-1$, $r_0= 0$ and there exist $k_1$, $k_2$ such that
$$
i_{r_{k_1}}+r_{k_1}-(\ell_1-1)=i_{\ell_1-1},
\quad
i_{r_{k_2}}+r_{k_2}=i_0.
$$
\end{itemize}

\begin{prop}\label{prop9-1}
After finitely many applications of the maps (\ref{eq9-12}), each element
is mapped to zero.
\end{prop}
\vskip2mm

\noindent
\pf
Let us prove the proposition by the induction on $\ell_1$.
For $\ell_1=1$, there is nothing to prove.

Suppose that $\ell_1>1$ and the statement holds until $\ell_1-1$.
Let us prove the $\ell_1$ case by the induction on 
$h=h(i_{\ell_1-1},\cdots,i_0)$.

If $h\leq \ell_1-2$, there is no $(i_{\ell_1-1},\cdots,i_0)$ satisfying
$i_{\ell_1-1}<\cdots<i_0$.

Suppose that $h\geq \ell_1-1$. 
We denote the value of $h$ of the right hand side of 
(\ref{eq9-12}) by $h'$. By Lemma \ref{lem9-2} $h'\leq h$.
If $h'<h$, the statement holds by the induction hypothesis.
Let us consider the case $h'=h$ which is divided into (i) to (iv) above.

We first consider the case (iv).

In this case $p_{\sigma_{k_1}}=\ell_1-1$, $p_{\sigma_{k_2}}=0$.
Therefore
\bea
&&
\sigma_{k_1}=k,
\quad
\sigma_{k_2}=0,
\quad
p_{k}=\ell_1-1,
\quad
p_0=0.
\non
\ena
We have
\bea
&&
i_{r_k}+r_k-p_{\sigma_k}
=i_{\ell_1-1}+\ell_1-1-p_{\sigma_k}
=i_{r_{k_1}}+r_{k_1}-p_{\sigma_k},
\non
\\
&&
i_{r_0}+r_0-p_{\sigma_0}
=i_0-p_{\sigma_0}
=i_{r_{k_2}}+r_{k_2}-p_{\sigma_0}.
\non
\ena
For the value of $\sigma_0$, $\sigma_k$ we have four cases.
\begin{itemize}
\item[(a).] $k_1=k$, $k_2=0$:
We have $\sigma_s\neq 0,k$ for any $1\leq s\leq k-1$.
\item[(b).] $k_1=k$, $k_2\neq 0$: We have $\sigma_0\neq 0,k$.
\item[(c).] $k_1\neq k$, $k_2=0$: We have $\sigma_k\neq 0,k$.
\item[(d).] $k_1\neq k$, $k_2\neq 0$: We have $\sigma_0,\sigma_k\neq 0,k$.
\end{itemize}

It follows that the map (\ref{eq9-12}) induces the map
\bea
&&
(i_{\ell_1-2},\cdots,i_1|n-r_{k-1},\cdots,n-r_1)
\non
\\
&&
\mapsto
\pm
(\cdots,i_{r_{k-1}}+r_{k-1}-p_{\sigma'_{k-1}},\cdots,
i_{r_{1}}+r_{1}-p_{\sigma'_{1}},\cdots|n-p_{k-1},\cdots,n-p_1),
\label{eq9-16}
\ena
where $\sigma'$ is the permutation of $(1,\cdots,k-1)$ specified 
according as the cases (a)-(d) above as
\begin{itemize}
\item[(a).] $\sigma_s'=\sigma_s$ for $1\leq s\leq k-1$,
\item[(b).] $\sigma_{k_2}'=\sigma_0$, $\sigma_s'=\sigma_s$ (otherwise),
\item[(c).] $\sigma_{k_1}'=\sigma_k$, $\sigma_s'=\sigma_s$ (otherwise),
\item[(d).] $\sigma_{k_1}'=\sigma_k$, $\sigma_{k_2}'=\sigma_0$, 
$\sigma_s'=\sigma_s$ (otherwise).
\end{itemize}

Notice that
\bea
&&
1\leq r_1<\cdots<r_{k-1}\leq \ell_1-2,
\quad
1\leq p_1<\cdots<p_{k-1}\leq \ell_1-2,
\non
\ena
and
\bea
&&
\text{LHS of (\ref{eq9-16})}=\pm(\text{RHS of (\ref{eq9-16})})
\non
\ena
is equivalent to
\bea
&&
\text{LHS of (\ref{eq9-12})}=\pm(\text{RHS of (\ref{eq9-12})}).
\non
\ena

Let us set
\bea
&&
I_{s-1}=i_s,
\quad
r_s-1=r'_{s-1},
\quad
p_s-1=p'_{s-1},
\quad
\sigma'_{s}-1=\tau_{s-1}.
\non
\ena
Then
\bea
&&
i_{r_s}+r_s-p_{\sigma'_s}=I_{r'_{s-1}}+r'_{s-1}-p'_{\tau_{s-1}},
\non
\\
&&
0\leq r'_1<\cdots<r'_{k-2}\leq \ell_1-3,
\quad
0\leq p'_1<\cdots<p'_{k-1}\leq \ell_1-3,
\non
\ena
and $\tau$ is the permutation of $(0,\cdots,k-2)$.
We have
\bea
&&
\text{LHS of (\ref{eq9-16})}=
(I_{\ell_1-3},\cdots,I_0|n-1-r'_{k-2},\cdots,n-1-r'_0),
\non
\\
&&
\text{RHS of (\ref{eq9-16})}=\pm
(\cdots,I_{r'_{k-2}}+r'_{k-2}-p'_{\tau_{k-2}},\cdots,
I_{r'_{0}}+r'_{0}-p'_{\tau_{0}},\cdots|
\non
\\
&&
\qquad\qquad\qquad\qquad
n-1-p'_{k-2},\cdots,n-1-p'_0),
\non
\ena
Thus the map (\ref{eq9-16}) can be considered as the map (\ref{eq9-12})
with the replacements
\bea
&&
\ell_1\ra \ell_1-2,
\quad
n\ra n-1,
\quad
k\ra k-2.
\non
\ena
Notice that, in the case (iv), the map beginning from 
the image of the map (\ref{eq9-12})
is again of the form (iv) if the image is not zero 
and the value of $h$ is invariant.
Therefore, by the induction hypothesis on $\ell_1$, if we apply the maps
(\ref{eq9-12}) of the type (iv) sufficiently many times to
an element we have zero. Thus the statement is proved in this case.

The cases (ii) and (iii) can be proved in a similar manner.
Things we should notify here are the followings.
In the case (ii) if the image is non-zero and the next map preserves
the value of $h$, the next map is of type (ii) or (iv).
Similarly for the case (iii) the next map is of type (iii) or (iv) under
the same condition.

Let us consider the case (i).
We divide the case into three.
\begin{itemize}
\item[(i-i).] $p_k\neq \ell_1-1$, $p_0=0$: The next map is in the case (ii)
for which the statement is already proved as above.

\item[(i-ii).] $p_k= \ell_1-1$, $p_0\neq 0$: The next map is in the case (iii)
for which the statement is already proved as above too.

\item[(i-iii).] $p_k\neq \ell_1-1$, $p_0\neq 0$:
\end{itemize}
In this case $k\leq \ell_1-3$ and the map (\ref{eq9-12}) induces the map
\bea
&&
(i_{\ell_1-2},\cdots,i_1|n-r_k,\cdots,n-r_0)
\non
\\
&&
\mapsto
(i_{\ell_1-2},\cdots,
i_{r_k}+r_k-p_{\sigma_k},\cdots,i_{r_0}+r_0-p_{\sigma_0},\cdots,
i_1|n-p_k,\cdots,n-p_0).
\label{eq9-17}
\ena
We set 
\bea
&&
I_{s-1}=i_s,
\quad
r_s-1=r'_s,
\quad
p_s-1=p'_s.
\non
\ena
Then
\bea
&&
\text{LHS of (\ref{eq9-17})}=(I_{\ell_1-3},\cdots,I_0|n-1-r'_k,
\cdots,n-1-r'_0),
\non
\\
&&
\text{RHS of (\ref{eq9-17})}=
(I_{\ell_1-3},\cdots,
I_{r'_{k}}+r'_k-p'_{\sigma_k},\cdots,I_{r'_{0}}+r'_0-p'_{\sigma_0},\cdots
|
\non
\\
&&
\qquad\qquad\qquad\quad
n-1-p'_k,\cdots,n-1-p'_0).
\non
\ena
Thus (\ref{eq9-17}) can be considered as the map (\ref{eq9-12}) with 
the modification
\bea
&&
\ell_1\ra \ell_1-2,
\quad
n\ra n-1.
\non
\ena
Therefore sufficiently many compositions of the maps of type (i-iii)
leads each element to zero by the induction hypothesis.

Let us summarize what is proved in the case (i).
To this end we introduce the case (v) which is (ii) or (iii) or (iv).
Then we have the diagram like
$$
\begin{array}{ccccccccccccc}
(i)&\ra&(v)&\ra&(v)&\ra&\cdots&\cdots&\cdots&\cdots&\cdots&\ra&0\\
{}&\sea&{}&{}&{}&{}&{}&{}&{}&{}&{}&{}\\
{}&{}&(i)&\ra&(v)&\ra&\cdots&\cdots&\cdots&\cdots&\cdots&\ra&0\\
{}&{}&{}&\sea&{}&{}&{}&{}&{}&{}&{}&{}&{}\\
{}&{}&{}&{}&\ddots&{}&{}&{}&{}&{}&{}&{}&{}\\
{}&{}&{}&{}&{}&\sea&{}&{}&{}&{}&{}&{}&{}\\
{}&{}&{}&{}&{}&{}&(i)&\ra&(v)&\ra&\cdots&\ra&0\\
{}&{}&{}&{}&{}&{}&{}&\sea&{}&{}&{}&{}&{}\\
{}&{}&{}&{}&{}&{}&{}&{}&0&{}&{}&{}&{},\\
\end{array}
$$
where the maps directed from NW to SE represent the composition of maps
of type (i-iii) and those from W to E represent the composition of maps of
type (v). This shows that in the case (i) the statement is proved.
Thus for $\ell_1$ the statement is proved by the induction on $h$.
Therefore the proposition is proved by the induction on $\ell_1$.
Q.E.D.
\vskip4mm

\noindent
Finally let us complete the proof of Theorem \ref{independence}.
For $k=0$, (\ref{eq9-10}) holds.
Then by Proposition \ref{prop9-1},
\bea
&&
\text{LHS of (\ref{eq9-10})}\equiv 0,
\label{eq9-18}
\ena
which proves (\ref{eq9-1}) for $k=0$.
If $k\geq 1$, (\ref{eq9-4}) holds by the induction hypothesis.
Consequently (\ref{eq9-10}) is valid. Again by Proposition \ref{prop9-1}
we have (\ref{eq9-18}) which proves (\ref{eq9-1}) for $k$.
Thus the proof of Theorem \ref{independence} is completed.
Q.E.D.

\section{Fermionic character formula}
For a graded vector space $W=\oplus_k W(k)$ such that $W(k)$ is 
finite dimensional, we define its character by
\bea
&&
\ch\, W=\sum_{k} q^k\,\dim W(k).
\non
\ena
We introduce a grading on $U_{2n,\ell}$ by $\deg_2$:
\bea
&&
U_{2n,\ell}=\oplus_{k=0}^\infty U_{2n,\ell}(k),
\ena
where $U_{2n,\ell}(k)$ is the subspace of elements with degree $k$.
By Proposition \ref{q-tetra-id}, Theorem \ref{conj-basis}, (\ref{ch2nl-1}) and
(\ref{ch2nl-2}) we have
\bea
&&
\ch\, U_{2n,\ell}=q^{n^2}\ch\, R_{2n}\cdot \ch_{2n,\ell}
=\frac{q^{n^2}}{[2n]_q!}
\qbc{2n}{\ell}.
\label{unl-char}
\ena

We consider the $R_{2n}$-module
\bea
&&
M_{2n,\ell}=\frac{U_{2n,\ell}}{U_{2n,\ell-1}\wedge w_1+U_{2n,\ell-2}\wedge \xi_1},
\label{Mnl}
\ena
here and in the following of this section we omit $(2n)$ of $v^{(2n)}_i$ etc.
The denominator of the right hand side of (\ref{Mnl}) is introduced because for any element $P$ in the denominator, 
$\psi_P$ vanishes identically due to Theorem \ref{tar} (ii), (\ref{w1=Xi1}) 
and (\ref{Xi2=xi1}).

Since $w_1$ and $\xi_1$ are homogeneous, $M_{2n,\ell}$ is also graded.
In this section we shall calculate the character of $M_{2n,\ell}$.
To this end we construct an $R_{2n}$ free resolution of $M_{2n,\ell}$.

Define the $R_{2n}$-linear map $\varphi_\ell$ by
\bea
&&
\varphi_\ell:\, 
U_{2n,\ell}\oplus U_{2n,\ell-1} \longrightarrow U_{2n,\ell+1}\oplus U_{2n,\ell},
\non
\\
&&
\varphi_\ell(a,b)=(a\wedge w_1+(-1)^\ell \, b\wedge \xi_1,\, b\wedge w_1),
\non
\ena
and the map $\psi_\ell$ by
\bea
&&
\psi_\ell:\, U_{2n,\ell}\oplus U_{2n,\ell-1} \longrightarrow U_{2n,\ell+1},
\non
\\
&&
\psi_\ell(a,b)=a\wedge w_1+(-1)^\ell\, b\wedge \xi_1,
\non
\ena
where we set $U_{n,0}=R_n$ and $U_{n,\ell}=0$ for $\ell<0$.
We denote $p_\ell$ the natural projection:
\bea
&&
p_\ell:\, U_{2n,\ell} \longrightarrow M_{2n,\ell}.
\non
\ena
The relation
\bea
&&
\varphi_\ell\varphi_{\ell-1}=0
\label{complex}
\ena
can be easily verified.

\begin{theorem}\label{resol}
The following sequence is exact for $0\leq \ell\leq n$:
\bea
&&
0\lar 
U_{2n,0}
\stackrel{\varphi_0}{\lar} 
U_{2n,1}\oplus U_{2n,0}
\stackrel{\varphi_1}{\lar}  
\cdots
\stackrel{\varphi_{\ell-2}}{\lar} 
U_{2n,\ell-1}\oplus U_{2n,\ell-2}
\stackrel{\psi_{\ell-1}}{\lar} 
U_{2n,\ell}
\stackrel{p_\ell}{\lar} 
M_{2n,\ell}
\lar
0.
\non
\ena
\end{theorem}
\vskip2mm

\noindent
\pf
We first prove 

\begin{lemma}\label{lem80}
The sequence
\bea
&&
U_{2n,m-1}\oplus U_{2n,m-2}
\stackrel{\varphi_{m-1}}{\lar} 
U_{2n,m}\oplus U_{2n,m-1}
\stackrel{\varphi_{m}}{\lar} 
U_{2n,m+1}\oplus U_{2n,m}
\label{eq-80}
\ena
is exact at the middle term for any $0\leq m\leq n-1$.
\end{lemma}
\vskip2mm

\noindent
\pf 
We have to prove
\bea
&&
\Ker\, \varphi_m=\Im \, \varphi_{m-1}.
\non
\ena

\begin{lemma}\label{lem81}
For any $0\leq m\leq n-1$, we have
\bea
&&
\ch\,(\Im \varphi_m)=\ch\, U_{2n,m}.
\non
\ena
\end{lemma}
\vskip2mm

\noindent
\pf 
We set
\bea
&&
U'_{2n,m+1}=
\oplus_{1\in J} R_{2n} v_I\wedge w_J\wedge \xi_K\subset U_{2n,m+1},
\non
\\
&&
U''_{2n,m-1}=
\oplus_{1\notin J} R_{2n} v_I\wedge w_J\wedge \xi_K\subset U_{2n,m-1},
\non
\ena
where, in both cases, the sum is taken for $I,J,K$ satisfying 
(\ref{indexI})-(\ref{indexK}). 
Define the map 
\bea
&&
\tilde{\varphi}_m:\, 
U'_{2n,m+1}\oplus U''_{2n,m-1}
\lar
\Im\, \varphi_m,
\non
\\
&&
\tilde{\varphi}_m(a,b)=(a+(-1)^mb\wedge \xi_1,\, b\wedge w_1).
\non
\ena
Notice the following three properties of the basis 
$Bas_\ell=\{v_I\wedge w_J\wedge \xi_K\}$ of $U_{2n,\ell}$ given in Theorem
\ref{conj-basis}.

The first is that, if $v_I\wedge w_J\wedge \xi_K$ is in $Bas_{\ell-1}$ 
and $1\notin J$, then $v_I\wedge w_1\wedge w_J\wedge \xi_K$ is in 
$Bas_{\ell}$. We denote the subset of $Bas_\ell$ consisting of such
elements by $Bas_\ell^{+}$.

The second is that, if $v_I\wedge w_J\wedge \xi_K$ is in $Bas_{\ell+1}$ 
and $J=(1,J')$, then $v_I\wedge w_{J'}\wedge \xi_K$ is in 
$Bas_{\ell}$. We denote the subset of $Bas_\ell$ consisting of such
elements by $Bas_\ell^{-}$.

The third property is 
\bea
&&
Bas_\ell=Bas_\ell^{+}\sqcup Bas_\ell^{-}.
\label{property-3}
\ena
It follows that $\tilde{\varphi}_m$ is bijective. 
Moreover $\tilde{\varphi}_m$ preserves the degree sinec 
$\deg_1\, w_1=\deg_1 \xi_1=0$.
Thus 
\bea
&&
\ch\,(\Im \,\varphi_m)=\ch\, U'_{2n,m+1}+\ch\, U''_{2n,m-1}=\ch\, U_{2n,m},
\non
\ena
where we use (\ref{property-3}) to derive the last equation.
Thus the lemma is proved.
Q.E.D.

We have the exact sequence:
\bea
&&
0\lar
\Ker\, \varphi_m
\lar
U_{2n,m}\oplus U_{2n,m-1}
\lar
\Im\, \varphi_{m}
\lar 0.
\ena
Then
\bea
&&
\ch\,(\Im\, \varphi_{m})=
\ch\, U_{2n,m}+\ch\, U_{2n,m-1}-\ch\,(\Ker \varphi_m).
\label{eq-81}
\ena
For two polynomials $f(q)$, $g(q)$ of $q$ we denote
\bea
&&
f(q)\geq g(q),
\non
\ena
if all the coefficients of $f(q)-g(q)$ are non-negative.
By (\ref{complex}) we have
\bea
&&
\Im \,\varphi_{m-1} \subset \Ker\, \varphi_m.
\non
\ena
Then, by (\ref{eq-81}) and Lemma \ref{lem81}, we have
\bea
\ch\,(\Im \,\varphi_m)
&\leq&
\ch\, U_{2n,m}+\ch\, U_{2n,m-1}-\ch\,(\Im \,\varphi_{m-1})
\non
\\
&=& \ch\, U_{2n,m}=\ch\,(\Im \,\varphi_m).
\non
\ena
Thus 
\bea
&&
\Im \,\varphi_{m-1}=\Ker\, \varphi_m,
\ena
which proves Lemma \ref{lem80}.
Q.E.D.

Let us prove 
\bea
&&
\Ker \,\psi_{\ell-1}=\Im \,\varphi_{\ell-2}.
\label{eq-83}
\ena
It is sufficient to prove 
\bea
&&
\Ker\, \psi_{\ell-1}\subset\Im \,\varphi_{\ell-2}.
\non
\ena
Suppose that
\bea
&&
\psi_{\ell-1}(a,b)=a\wedge w_1+(-1)^{\ell-1}b\wedge \xi_1=0,
\quad
a\in U_{2n,\ell-1},
\quad
b\in U_{2n,\ell-2}.
\non
\ena
By taking exterior product with $w_1$ we have
\bea
&&
(w_1\wedge b)\wedge \xi_1=0.
\non
\ena
Since the map
\bea
&&
\xi_1\wedge:\,
\wedge^m H^{(2n)}\lar \wedge^{m+2} H^{(2n)},
\quad
0\leq m\leq n-1,
\non
\ena
is injective as proved in \cite{Tar} using the representation theory of
$sl_2$,
\bea
&&
w_1\wedge b=0.
\non
\ena
Thus we have
\bea
&&
\varphi_{\ell-1}(a,b)=(0,0).
\non
\ena
Therefore, by Lemma \ref{lem80}, we have
\bea
&&
(a,b)\in \Im \,\varphi_{\ell-2}.
\non
\ena
Consequently (\ref{eq-83}) is proved.

Finally the equation
\bea
&&
\Ker\, p_\ell=\Im\, \psi_{\ell-1}
\non
\ena
is obvious by the definition of $M_{2n,\ell}$.
Thus Theorem \ref{resol} is proved.
Q.E.D.

\begin{cor}\label{cor81}
We have
\bea
&&
\ch\, M_{2n,\ell}=
\frac{q^{n^2}}{[2n]_q!}
\Big(
\qbc{2n}{\ell}-\qbc{2n}{\ell-1}
\Big).
\non
\ena
\end{cor}
\vskip2mm

\noindent
\pf
By Theorem \ref{resol}
\bea
\ch\, M_{2n,\ell}
&=&
\ch\, U_{2n,\ell}-
\sum_{r=0}^{\ell-1}(-1)^r
\big(
\ch\, U_{2n,\ell-1-r}+\ch\, U_{2n,\ell-2-r}
\big)
\non
\\
&=&
\ch\, U_{2n,\ell}-\ch\, U_{2n,\ell-1}.
\non
\ena
Then the corollary follows from (\ref{unl-char}).
Q.E.D.

For a non-negative integer $\lambda$ let us consider
the direct sum of $M_{2n,\ell}$ with $n-\ell$ being fixed:
\bea
&&
M^{(0)}_{2\lambda}:=\oplus_{n-\ell=\lambda}M_{2n,\ell}.
\non
\ena
Then, by Corollary \ref{cor81}, we get
\bea
&&
\ch\, M^{(0)}_{2\lambda}=
\sum_{n-\ell=\lambda}\frac{q^{n^2}}{[2n]_q!}
\Big(
\qbc{2n}{\ell}-\qbc{2n}{\ell-1}
\Big).
\label{ch-m2lambda}
\ena

Let $\Lambda_i$, $i=0,1$ be the fundamental weights of $\widehat{sl_2}$,
$V(\Lambda_i)$ be the level one irreducible highest weight representation
with the highest weight $\Lambda_i$,
\bea
&&
e_i, 
\quad
f_i, 
\quad
h_i,
\quad
i=0,1
\non
\ena
the Chevalley generators of $\widehat{sl_2}$ and $d$ the scaling element 
\cite{Kac}.
Consider the subspace of $V(\Lambda_i)$ consisting of $sl_2$ highest
weight vectors with the $sl_2$ weight $\mu$:
\bea
&&
V(\Lambda_i|\mu)=
\{v\in V(\Lambda_i)\, |\, e_1 v=0,\,\, h_1v=\mu v\,\}.
\non
\ena
Then
\bea
&&
\ch\, M^{(0)}_{2\lambda}=\tr_{V(\Lambda_0|2\lambda)}(q^{-d}),
\non
\ena
which is the branching function of $2\lambda+1$-dimensional irreducible
representation of $sl_2$ in the representation $V(\Lambda_0)$ of
$\widehat{sl_2}$. 
In fact if we multiply (\ref{ch-m2lambda}) by the character $\chi_{\lambda}(z)$
of the $2\lambda+1$-dimensional irreducible representation of $sl_2$ we get
\bea
&&
\sum_{\lambda=0}^\infty
\chi_{\lambda}(z)\, \ch\, M^{(0)}_{2\lambda}
=
\sum_{n=0}^\infty\sum_{\ell=0}^{2n}
\frac{z^{2n-2\ell}q^{n^2}}{[\ell]_q!\,[2n-\ell]_q! },
\non
\ena
which is the fermionic form of the character of $V(\Lambda_0)$
presented in \cite{KMM,Me}. 

The space $V(\Lambda_0|2\lambda)$ becomes the irreducible representation 
of the Virasoro algebra with $c=1$ and $h=\lambda^2$ \cite{Kac}.
It follows that \cite{Roch}
\bea
&&
\tr_{V(\Lambda_0|2\lambda)}(q^{-d})=
\frac{q^{\lambda^2}(1-q^{2\lambda+1})}{(q:q)_\infty},
\quad
(z:q)_\infty=\prod_{j=0}^\infty(1-zq^j).
\non
\ena

\section{The case of odd number of variables}
In this section we study the case where the number of $x_j$'s is odd.
The structure is similar to the even case. But there are some differences
in the construction of a basis. Thus we consider this case separately.
We omit the proofs of the statements if they are quite similar to the even case.

Let us define $P^{(2n+1)}_{r,s}$, $1\leq r\leq n+1$, $s\in \mathbb{Z}$
by the same recursion relation as (\ref{def-P}):
\bea
P^{(2n+1)}_{1,s}&=&e^{(2n+1)}_{2s-1},
\non
\\
P^{(2n+1)}_{r,s}&=&P^{(2n+1)}_{r-1,s+1}-e^{(2n+1)}_{2s}P^{(2n+1)}_{r-1,1}
\quad \hbox{for $r\geq 2$}.
\label{def-o-P}
\ena
Then $P^{(2n+1)}_{r,s}=0$ for $s\leq 0$.

\begin{prop}\label{op-1}
\begin{itemize}
\item[(i).] 
$P^{(2n+1)}_{r,s}=0$ for $r\geq 2$ and $s\geq n+1$.

\item[(ii).] 
$\overline{P^{(2n+1)}_{r,s}}=P^{(2n-1)}_{r,s}-x^2P^{(2n-1)}_{r,s-1}$.
\end{itemize}
\end{prop}
\vskip3mm

Let us define $v^{(2n+1)}_i$, $w^{(2n+1)}_i$, $0\leq i\leq n$ by
\bea
&&
v^{(2n+1)}_0=\sum_{j=0}^n e^{(2n+1)}_{2j}X^{2j},
\quad
v^{(2n+1)}_i=\sum_{j=1}^{n+1} P_{i,j}^{(2n+1)}X^{2(j-1)}
\,\,(i\neq 0),
\quad
w^{(2n+1)}_i=\sum_{j=1}^{n+1} P_{i+1,j}^{(2n+1)}X^{2j-1}.
\non
\ena
By Proposition \ref{op-1}
\bea
&&
\rho_{\pm}(v^{(2n+1)}_i)=\rho_{\pm}(w^{(2n+1)}_i)=0,
\quad
0\leq i\leq n.
\non
\ena
Notice that 
\bea
&&
w^{(2n+1)}_0\notin H^{(2n+1)},
\quad
v^{(2n+1)}_i, w^{(2n+1)}_j\in H^{(2n+1)},
\quad j\neq 0.
\non
\ena
We set, for $1\leq k\leq n$,
\bea
2\xi_k^{(2n+1)}(X_1,X_2)&=&
\frac{X_1-X_2}{X_1+X_2}
\Big(
v_0^{(2n+1)}(X_1)w^{(2n+1)}_{k-1}(X_2)+
v_0^{(2n+1)}(X_2)w^{(2n+1)}_{k-1}(X_1)
\Big)
\non
\\
&&
+v_0^{(2n+1)}(X_1)w^{(2n+1)}_{k-1}(X_2)-v_0^{(2n+1)}(X_2)w^{(2n+1)}_{k-1}(X_1).
\label{def-0-xi}
\ena
Then 
\bea
&&
\xi_k^{(2n+1)}\in H^{(2n+1)},
\quad
\rho_\pm(\xi_k^{(2n+1)})=0,
\non
\ena
which means
$$
\xi_k^{(2n+1)}\in U_{2n+1,2}. 
$$
We have the relation
\bea
&&
\Xi^{(2n+1)}_1=v_0^{(2n+1)},
\quad
\Xi^{(2n+1)}_2=\xi_1^{(2n+1)}.
\non
\ena

We set 
\bea
&&
Bas_\ell^{o-}=\{v^{(2n+1)}_I\wedge w^{(2n+1)}_J\wedge \xi^{(2n+1)}_K\},
\quad
Bas_\ell^{o+}=v_0^{(2n+1)}\wedge Bas_{\ell-1}^{o-},
\non
\ena
where $I,J,K$ satisfy (\ref{indexI})-(\ref{sum=l}) and
define
\bea
&&
Bas_\ell^{o}=Bas_\ell^{o+}\sqcup Bas_\ell^{o-}.
\non
\ena
It follows from the even case that the cardinality of $Bas_\ell^{o}$ is
\bea
&&
\sharp Bas_\ell^{o}
=
\sharp Bas_{\ell-1}^{o-}+\sharp Bas_\ell^{o-}
=
\bc{2n}{\ell-1}+\bc{2n}{\ell}
=
\bc{2n+1}{\ell}.
\non
\ena
The degrees of elements are given by
\bea
&&
\deg_1\, v^{(2n+1)}_0=0,
\quad
\deg_1\, v^{(2n+1)}_i=2i-1\,\,
(i\neq 0),
\quad
\deg_1\, w^{(2n+1)}_j=2j,
\quad
\deg_1\, \xi^{(2n+1)}_k=2k-2.
\non
\ena
The structure of $U_{2n+1,\ell}$ is given by the following theorem.

\begin{theorem}\label{oth-1}
\begin{itemize}
\item[(i).] 
The module $U_{2n+1,\ell}$ is a free $R_{2n+1}$-module of rank
$\binom{2n+1}{\ell}$ with $Bas_\ell^{o}$ as a basis.
\item[(ii).] 
The character $\ch_{2n+1,\ell}$ of $U_{2n+1,\ell}$ with respect to 
$\deg_1$is given by
\bea
\ch_{2n+1,\ell}=\qbc{2n+1}{\ell}.
\non
\ena
\end{itemize}
\end{theorem}

Let us first prove (ii) assuming (i).
Let 
\bea
&&
b_{n,\ell}^{\pm}=\sum_{\ell_1+\ell_2+2\ell_3=\ell}
q^{\ell_1^2+\ell_2(\ell_2\pm 1)}
\qqtetc{n}{\ell_1}{\ell_2}{\ell_3}.
\non
\ena
By Proposition \ref{q-tetra-id} we have
\bea
&&
b_{n,\ell}^{-}=\qbc{2n}{\ell}.
\non
\ena
We easily have
\bea
&&
\ch_{2n+1,\ell}=b_{n,\ell-1}^{+}+b_{n,\ell}^{+},
\ena
where in the right hand side the first and the second terms
are the characters of $Bas_\ell^{o+}$ and $Bas_\ell^{o-}$ respectively.
Notice a similar identity to (\ref{qtetra-id2}):
\bea
&&
[\ell_2]_{q^2}
\qqtetc{n}{\ell_1}{\ell_2}{\ell_3}
=
[\ell_3+1]_{q^2}
\qqtetc{n}{\ell_1}{\ell_2-1}{\ell_3+1},
\non
\ena
which implies
\bea
q^{2\ell_2}\qqtetc{n}{\ell_1}{\ell_2}{\ell_3}
&=&
\qqtetc{n}{\ell_1}{\ell_2}{\ell_3}
-
\qqtetc{n}{\ell_1}{\ell_2-1}{\ell_3+1}
\non
\\
&&
+
q^{2(\ell_3+1)}
\qqtetc{n}{\ell_1}{\ell_2-1}{\ell_3+1}.
\non
\ena
Using this relation we rewrite $b_{n,\ell}^{+}$ as
\bea
b_{n,\ell}^{+}&=&
\sum_{\ell_1+\ell_2+2\ell_3=\ell}
q^{\ell_1^2+\ell_2(\ell_2-1)+2\ell_2}
\qqtetc{n}{\ell_1}{\ell_2}{\ell_3}
\non
\\
&=&
\qbc{2n}{\ell}
\non
\\
&&
-\sum_{\ell_1+\ell_2+2\ell_3=\ell+1, \ell_3\geq 1}
q^{\ell_1^2+\ell_2(\ell_2+1)}
\qqtetc{n}{\ell_1}{\ell_2}{\ell_3}
\label{oeq-1}
\\
&&
+\sum_{\ell_1+\ell_2+2\ell_3=\ell+1, \ell_3\geq 1}
q^{\ell_1^2+\ell_2(\ell_2+1)+2\ell_3}
\qqtetc{n}{\ell_1}{\ell_2}{\ell_3}.
\label{oeq-2}
\ena
We easily have
\bea
&&
(\ref{oeq-1})+(\ref{oeq-2})=q^{\ell+1}\qbc{2n}{\ell+1}-b_{n,\ell+1}^{+}.
\non
\ena
The desired result follows from (\ref{binom-rec2}).
Q.E.D

The proof of (i) consists of two parts, the linear independence and 
the generation of the space $U_{2n+1,\ell}$. Let us first prove the linear
independence. It is proved by the specialization argument as in the even case
using the result of the even case.

Consider the specialization of variables:
\bea
&&
e_{2n+1}^{(2n+1)}=-e_{2n}^{(2n+1)}=1,
\quad
e_{j}^{(2n+1)}=0,
\quad
j\neq 2n,2n+1.
\label{o-special}
\ena

\begin{prop}\label{op-2}
At (\ref{o-special}) we have
\begin{itemize}
\item[(i). ]
$$
P^{(2n+1)}_{r,s}=
\left\{
\begin{array}{rl}
1,&\quad \text{$s=n+2-r$}\\
0,&\quad \text{otherwise},
\end{array}
\right.
$$

\item[(ii).]
$$
v^{(2n+1)}_0=1-X^{2n}, 
\quad
v^{(2n+1)}_i=X^{2(n+1-i)}, 
\quad
i\neq 0,
\quad
w^{(2n+1)}_j=X^{2n+1-2j},
$$

\item[(iii).] 
$$
\xi^{(2n+1)}_k=-\sum_{r=1}^n v^{(2n+1)}_{k-r}\wedge w^{(2n+1)}_r,
$$
where the index $k-r$ of $v^{(2n+1)}_{k-r}$ is read by modulo $n$ in the 
representative $\{1,2,...,n\}$.
\end{itemize}
\end{prop}

At (\ref{o-special}) we define $\gamma$, $\alpha_i$, $\beta_i$, $\omega_i$, 
$1\leq i\leq n$ by
\bea
&&
\gamma=v^{(2n+1)}_0,
\quad
\alpha_1= v^{(2n+1)}_n,
\quad
\alpha_i= v^{(2n+1)}_{n+1-i},
\quad
i\geq 2,
\non
\\
&&
\beta_i=w_i^{(2n+1)},
\quad
\omega_i=\sum_{r=1}^n\alpha_{r-i+1}\wedge \beta_r,
\non
\ena
where the indices of $\alpha_i$, $\beta_j$ are considered by modulo $n$ in the representative $\{1,2,\cdots,n\}$. 
By Proposition \ref{op-2}, $\{\alpha_i,\beta_j,\gamma\}$ are linearly independent.

Corresponding to $Bas_\ell^{o-}$ we define
\bea
&&
SBas_\ell^{o-}=
\{\alpha_I^{(2n+1)}\wedge \beta_J^{(2n+1)}\wedge \omega_K^{(2n+1)}\},
\non
\ena
where $I,J,K$ satisfy (\ref{indexI})-(\ref{sum=l}). 
By Theorem \ref{equiv-cond} and \ref{independence},
elements of $SBas_\ell^{o-}$ are linearly independent.
Thus elements of
\bea
&&
\gamma\wedge SBas_{\ell-1}^{o-} \sqcup SBas_\ell^{o-}
\non
\ena
are linearly independent. Consequently elements of $Bas_\ell^o$ are linearly independent over $R_{2n+1}$. 
\vskip2mm

Next let us prove that elements of $Bas_\ell^o$ generate $U_{2n+1,\ell}$ over $R_{2n+1}$. 
Number the elements in $Bas_\ell^o$ from $1$ to $\binom{2n+1}{\ell}$ and
name the $r$-th element $Q_r$. Expand $Q_r$ as in (\ref{qr}) and define the 
$\binom{2n+1}{\ell}$ by $\binom{2n+1}{\ell}$ matrix $X^{(2n+1,\ell)}$ by the similar formula to (\ref{matx}).
The assertion follows from

\begin{prop}\label{op-4}
We have
\bea
&&
\det\,X^{(2n+1,\ell)}=c\cdot (\Delta^{+}_{2n+1})^{\binom{2n}{\ell-1}+\binom{2n-1}{\ell-1}},
\non
\ena
for some non-zero constant $c$.
\end{prop}

Let 
\bea
&&
d^o_{I,J,K}=\deg_1\, v^{(2n+1)}_I\wedge w^{(2n+1)}_J \wedge \xi^{(2n+1)}_K,
\ena
and
\bea
&&
d^{(2n+1,\ell)}=
\sum_{I,J,K}\, d^o_{I,J,K}
+\sum_{0\leq s_1\leq\cdots\leq s_\ell\leq 2n}(s_1+\cdots+s_\ell),
\non
\ena
where $(I,J,K)$ runs index sets such that 
$ v^{(2n+1)}_I\wedge w^{(2n+1)}_J \wedge \xi^{(2n+1)}_K\in Bas_\ell^{o}$.
Then Proposition \ref{op-4} follows from 

\begin{lemma}\label{olem-1}
\bea
&&
d^{(2n+1,\ell)}=
\binom{2n+1}{2}\Biggl(\binom{2n}{\ell-1}+\binom{2n-1}{\ell-1}\Biggr).
\non
\ena
\end{lemma}
\vskip2mm

The proof of this lemma is similar to that of Lemma \ref{lem-d2nl}.

In this way (i) of Theorem \ref{oth-1} is proved.
Q.E.D.
\vskip3mm

In the following we omit $(2n+1)$ of $v^{(2n+1)}_i$ etc.

Let
\bea
&&
M_{2n+1,\ell}=\frac{U_{2n+1,\ell}}{U_{2n+1,\ell-1}\wedge v_0+U_{2n+1,\ell-2}\wedge \xi_1}.
\non
\ena
We define the maps $\varphi_\ell^o$, $\psi_\ell^o$, $p_\ell^o$ replacing $U_{2n,\ell}$ by $U_{2n+1,\ell}$, $w_1$ by $v_0$ and $M_{2n,\ell}$ by $M_{2n+1,\ell}$ in the definition of $\varphi_\ell$, $\psi_\ell$ and $p_\ell$.
We introduce a grading on $M_{2n+1,\ell}$ by $\deg_2$.

Then

\begin{theorem}\label{oth-2}

\begin{itemize}
\item[(i).] The following sequence is exact for $0\leq \ell\leq n$:
\bea
0\lar 
U_{2n+1,0}
\stackrel{\varphi_0^o}{\lar} 
U_{2n+1,1}\oplus U_{2n+1,0}
\stackrel{\varphi_1^o}{\lar}  
\cdots
&\stackrel{\varphi_{\ell-2}^o}{\lar}&
U_{2n+1,\ell-1}\oplus U_{2n+1,\ell-2}
\non
\\
&\quad&
\stackrel{\psi_{\ell-1}^o}{\lar}
U_{2n+1,\ell}
\stackrel{p_\ell^o}{\lar} 
M_{2n+1,\ell}
\lar
0.
\non
\ena

\item[(ii).] We have
\bea
&&
\ch\, M_{2n+1,\ell}=
\frac{q^{\frac{1}{4}(2n+1)^2}}{[2n+1]_q!}
\Big(
\qbc{2n+1}{\ell}-\qbc{2n+1}{\ell-1}
\Big).
\non
\ena
\end{itemize}
\end{theorem}

Similarly to the even case if we set
\bea
&&
M^{(1)}_{2\lambda+1}=\oplus_{n-\ell=\lambda}M_{2n+1,\ell},
\non
\ena
we have 
\bea
&&
\ch\, M^{(1)}_{2\lambda+1}=
\sum_{n-\ell=\lambda}
\frac{q^{\frac{1}{4}(2n+1)^2}}{[2n+1]_q!}
\Big(
\qbc{2n+1}{\ell}-\qbc{2n+1}{\ell-1}
\Big)
\non
\ena
which is equal to
\bea
&&
\tr_{V(\Lambda_1|2\lambda+1)}(q^{-d+\frac{1}{4}})=
\frac{q^{\frac{1}{4}(2\lambda+1)^2}(1-q^{2\lambda+2})}{(q:q)_\infty}.
\non
\ena

\appendix
\section{Proof of Proposition \ref{prop2-2}}
Let us write the integrand of $I_M(P)$ by $J_M(P)$:
\bea
&&
I_M(P)=\int_{C^\ell}\prod_{a=1}^\ell d\alpha_a J_M(P).
\non
\ena
We simply write $\Res$ instead of writing 
\bea
&&
\Res_{\beta_n=\beta_{n-1}+\pi i}.
\non
\ena
We divide the case into four:
(I). $n-1,n\notin M$, (II). $n-1\in M$, $n\notin M$,
(III) $n-1\notin M$, $n\in M$, (IV). $n-1\in M$, $n\in M$.
As in the proof of Proposition 3 in \cite{NT} one can calculate the residue
in the following way according as the four cases.

\noindent
(I). 
\bea
&&
\frac{1}{2\pi i}\Res\, I_M(P)=
Res\sum_{a=1}^\ell \int_{C^{\ell-1}}
(\Res_{\alpha_a=\beta_{n-1}-\pi i}J_M(P))\prod_{b\neq a} d\alpha_b.
\non
\ena

\noindent
(II).
\bea
\frac{1}{2\pi i}\Res\, I_M(P)&=&
Res\sum_{a=1}^\ell \int_{C^{\ell-1}}
(\Res_{\alpha_a=\beta_{n-1}-\pi i}J_M(P))\prod_{b\neq a} d\alpha_b
\non
\\
&&-
\Res \int_{C^{\ell-1}}
(\Res_{\alpha_\ell=\beta_{n-1}}J_M(P))\prod_{b=1}^{\ell-1} d\alpha_b.
\non
\ena

\noindent
(III).
\bea
\frac{1}{2\pi i}\Res\, I_M(P)&=&
Res\sum_{a=1}^\ell \int_{C^{\ell-1}}
(\Res_{\alpha_a=\beta_{n-1}-\pi i}J_M(P))\prod_{b\neq a} d\alpha_b
\non
\\
&&-
\Res \int_{C^{\ell-1}}
(\Res_{\alpha_\ell=\beta_{n-1}}J_M(P))\prod_{b=1}^{\ell-1} d\alpha_b
\non
\\
&&
+Res\sum_{a=1}^\ell \int_{C^{\ell-1}}
(\Res_{\alpha_a=\beta_{n-1}+\pi i}J_M(P))\prod_{b\neq a} d\alpha_b.
\non
\ena

\noindent
(IV).
\bea
\frac{1}{2\pi i}\Res\, I_M(P)&=&
Res\sum_{a=1}^\ell \int_{C^{\ell-1}}
(\Res_{\alpha_a=\beta_{n-1}-\pi i}J_M(P))\prod_{b\neq a} d\alpha_b
\non
\\
&&-
\Res \sum_{a=\ell-1,\ell}\int_{C^{\ell-1}}
(\Res_{\alpha_a=\beta_{n-1}}J_M(P))\prod_{b\neq a} d\alpha_b
\non
\\
&&
+Res\sum_{a=1}^\ell \int_{C^{\ell-1}}
(\Res_{\alpha_a=\beta_{n-1}+\pi i}J_M(P))\prod_{b\neq a} d\alpha_b.
\non
\ena
\vskip2mm

The theorem can be proved by the calculations using these formulae.

\end{document}